\documentclass[12pt]{amsart}
\usepackage{amsmath,amsthm,amsfonts,amssymb,enumitem,xcolor,verbatim}
\begin{document} 
\newcommand{\eqdef}{\overset{\mathrm{def}}{=\joinrel=}}
\newcommand{\C}{{\mathbb C}}
\newcommand{\dfe}{\overset{\mathrm{def}}{=\joinrel=}}
\newcommand{\N}{{\mathbb N}}
\renewcommand{\P}{{\mathbb P}}
\newcommand{\R}{{\mathbb R}}
\newcommand{\Aut}{\mathop{Aut}}
\newcommand{\Sing}{\mathop{Sing}}
\newtheorem{lemma}{Lemma}[section]
\newtheorem{definition}[lemma]{Definition}
\newtheorem{claim}[lemma]{Claim}
\newtheorem{corollary}[lemma]{Corollary}
\newtheorem*{conjecture}{Conjecture}
\newtheorem*{SpecAss}{Special Assumptions}
\newtheorem{example}{Example}
\newtheorem*{remark}{Remark}
\newtheorem*{observation}{Observation}
\newtheorem*{fact}{Fact}
\newtheorem*{remarks}{Remarks}
\newtheorem{proposition}[lemma]{Proposition}
\newtheorem{theorem}[lemma]{Theorem}
\newcommand{\cH}{{\mathcal H}}
\newcommand{\cTX}{{\mathcal T}_X}
\newcommand{\chH}{\widehat{\cH}}
\newcommand{\vv}{\alpha_i^*v}
\newcommand{\ww}{w_{ij_0}}

\numberwithin{equation}{section}
\def\labelenumi{\rm(\roman{enumi})}
\keywords{Tame sets, Danielewski surfaces, discrete subsets of complex
manifolds}
\title{%
  Rosay-Rudin-Spaces, Tame Sets and Danielewski Surfaces
}
\author {J\"org Winkelmann}
\begin{abstract}
  We introduce the notion of an ``$RR$-space''.
  ($RR$ stands for Rosay and Rudin.)
  These spaces essentially
  share the properties of tame subsets known for $\C^ n$.

  The class of $RR$-spaces contains character-free complex linear
  algebraic groups as well as Danielewski surfaces.
\end{abstract}
%
\address{%
J\"org Winkelmann \\
Lehrstuhl Analysis II \\
Mathematisches Institut \\
Ruhr-Universit\"at Bochum\\
44780 Bochum \\
Germany
}
\email{joerg.winkelmann@rub.de
}
%
\maketitle

\section{Introduction}

In their seminal paper \cite{RR} Rosay and Rudin proved remarkable
results on discrete subsets in $\C^n$. They called a discrete
subset $D$ of $\C^n$ {\em tame} if there is a holomorphic automorphism
$\phi$ of $\C^n$ with
\[
\phi(D)=\N\times\{0\}^{n-1}.
\]
There exist both tame and non-tame discrete subsets in $\C^n$.
Every injective self-map of a tame discrete subset extends to a
holomorphic automorphism of $\C^n$.

To extend these and other results on $\C^n$ to other manifolds,
{\em tameness} notions for other complex manifolds were
introduced.
In \cite{JW-TAME1} we introduced the notion of
(weakly) tame, while Andrist and Ugolini proposed a notion
of strongly tame(\cite{AU}).

For $\C^n$ all these tameness notions agree, but not for arbitrary
complex manifolds.

In this paper we propose a set of axioms for ``{\em RR-spaces''}
such that for these $RR$-spaces all the tameness notions agree and
we have (mostly) the same results as for $\C^n$.
We summarize the properties of $RR$-spaces in Theorem~\ref{summary}.

Besides $\C^n$, the class of $RR$-spaces contains all
character-free complex linear algebraic groups (basically this
is the main content of \cite{JW-LAG}), and
the class of Danielewski surfaces.

Danielewski surface (also called Gizatullin surfaces of rank $1$)
are a certain class of Stein surfaces which admit many automorphisms
(the action on the surface is transitive, and the automorphism group
is infinite-dimensional), but no (finite-dimensional) complex Lie group acts
transitively. Hence different methods are required to study
tame discrete subsets in these surfaces.

\section{RR-Spaces}

\subsection{Tameness}
Here we recall the different tameness notions.

Rosay and Rudin called a discrete subset $D$ of $\C^n$ {\em tame}
if there exists an holomorphic automorphism $\phi$ of $\C^n$
with $\phi(D)=\N\times\{0\}^{n-1}$.

In \cite{JW-TAME1} we called an infinite  discrete subset $D$ in
a complex manifold $X$   {\em (weakly) tame}
if for every  exhaustion function $\rho:X\to\R^+$
and every map $\zeta:D\to\R^+$ there exists a holomorphic automorphism $\phi$
of $X$ such that $\rho(\phi(x))\ge\zeta(x)$.

Following \cite{AU}, an infinite discrete subset $D$ of a complex manifold $X$
is called {\em (strongly) tame} if for every injective map $f:D\to D$
there exists a holomorphic automorphism $\phi$ of $X$ with $\phi|_D=f$.

For $\C^n$, all three notions are equivalent.
``Strongly tame'' evidently implies ``weakly tame''.
However, there are complex manifolds, e.g.,
\[
X=\{(z,w)\in\C^2:|z|<1\}
\]
with infinite discrete subsets which are weakly tame, but
not strongly tame. (\cite{JW-TAME1},\S5).

\subsection{Axioms}

Let us state the axioms for an ``$RR$-space''.

\begin{definition}
  Let $X$ be a complex manifold, ${\mathcal D}$ the set of all
  infinite discrete subsets of $X$.

  $X$ is a {\em $RR$-space} if there is a subset
  $\mathcal T\subset\mathcal D$ satisfying the conditions below.
  (Here a discrete subset $D\subset X$ is called
  ``(axiomatically) tame'' if $D\in\mathcal T$).

\begin{enumerate}[label=(RR\arabic*)]
\item Given a tame discrete set $D$ and a tame discrete set $D'$, every
injective map $f$ from $D$ to $D'$ extends to a biholomorphic self-map 
$\phi$ of $X$.
\item
There is an exhaustion function $\rho:X\to\R^+$ and an unbounded sequence of
positive real numbers $R_n$ such that every infinite subset $D\subset X$
with
\[
\#\{x\in D:\rho(x)\le R_n\}\le n \ \ \forall n\in\N
\]
must be a tame discrete subset.
\item
  If $\phi$ is a holomorphic automorphism of $X$ and $D\in \mathcal T$,
  then $\phi(D)\in\mathcal T$.
\item
  There is a non-zero complete vector field $v$ on $X$
  and a non-constant holomorphic function $f$ such that $v(f)=0$.
\end{enumerate}
\end{definition}

\begin{remark}
  Condition $(RR2)$ is equivalent to:
  
  {\em For every exhaustion function $\rho:X\to\R^+$
    there exists  unbounded sequence of
positive real numbers $R_n$ such that every subset $D\subset X$
with
\[
\#\{x\in D:\rho(x)\le R_n\}\le n \}\ \ \forall n\in\N
\]
must be a tame discrete subset.}
  See \cite{JW-TAME1}.
\end{remark}

\subsection{Consequences from the axioms}
Here we list some consequences of the axioms, i.e.,
properties satisfied by all $RR$-spaces.

\begin{theorem}\label{summary}
  Every $RR$-space  $X$ satisfies the following statements:
  \begin{enumerate}
  \item (Prop.~\ref{prop1})
    For discrete subsets $D$ of $X$ all tameness notions discussed
    above are equivalent, i.e., $D$ satisfies one of the following properties
    if and only if it satisfies all of these properties:
    \begin{itemize}
    \item
      $D\in \mathcal T$,
    \item
      $D$ is weakly tame.
    \item
      $D$ is strongly tame.
    \end{itemize}
  \item (Prop.~\ref{stein})
   $X$ is a Stein manifold.
  \item
    (Prop.~\ref{oka})
  Let $D$ be a discrete subset of $X$ which is
  empty, finite or tame.

  Then $X^*=X\setminus D$ is an Oka manifold
  in the sense of \cite{Forstneric}.
  \item
    (Prop.~\ref{not-loc-cpt})
    $\Aut(X)$ is not locally compact
    (for the compact-open topology).
  \item
    (Prop.~\ref{no-equ-rel})
    If $\sim$ is an equivalence relation on $X$ which is invariant under
    all $\phi\in\Aut^ 0(X)$
    (where $\Aut^0(X)$ denotes the connected component of the identity
    of the group of holomorphic automorphisms $\Aut(X)$
    endowed with the compact-open topology),
    then $\sim$ is trivial, i.e., either
    $x\sim y\ \forall x,y\in X$ or $\forall x,y\in X: x\sim y\iff x=y$.
  \item
    (Prop.~\ref{no-metric})
    The Kobayashi and Caratheodory pseudometrics on $X$ are
    trivial.
  \item (Corollary~\ref{no-bd-psh})
    Every bounded plurisubharmonic function on $X$ is constant.
  \item
(Prop.~\ref{t-span-vf})
    Let $D$ be a tame discrete set or empty or finite and $X^*=X\setminus D$.
    
    $\Aut^0(X^*)$  acts $N$-transitively on $X^*$
    for every $N\in\N$.
  \item
    (Prop.~\ref{non-tame})
    There exists an infinite discrete subset of $X$ which is not tame.
  \item
    (Corollary~\ref{cor-sub-tame})
    Every finite subset of $X$ is contained in a tame discrete subset.
  \item
    (Corollary~\ref{sym-diff})
    Let $D$ and $D'$ be discrete subsets of $X$ such that
    the ``symmetric difference'' $(D\setminus D')\cup(D'\setminus D)$ is
    finite.

    Then $D$ is tame if and only if $D'$ is tame.
  \item
    (Prop.~\ref{prop3})
    Every infinite subset of a tame discrete subset is itself tame.
  \item
    (Prop.~\ref{unbounded-contains-tame})
    Every unbounded subset of $X$ contains a tame discrete subset.
  \item
    (Prop.~\ref{avoid})
    For every tame discrete subset $D$ there are
    non-degenerate holomorphic maps from both $X$ itself
    and $\C^n$ ($n=\dim(X)$) to $X$  which avoid $D$.
  \item
    (Prop.~\ref{countable})
    Let $S$ and $S'$ be dense countable subsets of $X$. Then there
    exists an automorphism $\phi$ of $X$ with $\phi(S)=S'$.
\end{enumerate}
\end{theorem}

\subsection{Examples of $RR$-spaces}

\begin{theorem}
The following classes of complex manifolds satisfy all
axioms of an $RR$-space:

\begin{enumerate}
\item
  $\C^n$. This is, in our language, the main result of
  Rosay and Rudin in \cite{RR}.
\item
  Character-free complex linear algebraic groups, i.e., algebraic subgroups
  of $GL_n(\C)$ which do not admit a
  non-constant morphism of algebraic groups
  to the multiplicative group $\C^*$.
  (\cite{JW-LAG}, see Theorem~\ref{LAG-is-RR}).
  For example, this includes all semisimple complex Lie groups
  like $SL_n(\C)$.
\item
  Danielewski surfaces, i.e., surfaces which may be described as
  \[
  X=\{(x,y,z)\in\C^3: xy=P(z)\}
  \]
  for a polynomial $P\in\C[t]$ without multiple zeroes.
  (Theorem~\ref{dan-is-RR}).
\end{enumerate}
\end{theorem}

\subsection{Additional properties}

For some $RR$-spaces we are able to deduce some additional results.
\begin{theorem}
  Let $X$ be a character-free complex linear algebraic group (CFLAG) or
  a Danielewski surface.

  Then the following assertions hold:

  \begin{enumerate}
  \item
  Every infinite discrete subset is the union of two tame
  infinite discrete subsets.
  (CFLAG:
  \cite{JW-LAG}, Theorem~17;
  Danielewski surfaces: Corollary~\ref{dan-union-tame}).
 \item
   An infinite discrete subset $D$ of $X$ is tame if and only if
   there exists an algebraic curve $C\subset X$ and a holomorphic
   automorphism $\phi$ of $X$ such that $\phi(D)\subset C\simeq\C$.
   (Danielewski surfaces: Proposition~\ref{flex-tame}; CFLAG: \cite{JW-LAG}, Theorem~1)
  \end{enumerate}
\end{theorem}

We believe that both statements should be true for every $RR$-space $X$
($(ii)$ if $X$ is in addition an affine algebraic variety)
but are currently unable to deduce them from the $RR$-space axioms.

\subsection{Conjectures}

\begin{conjecture}
  \begin{enumerate}
  \item
    Every smooth flexible complex affine variety is an
    $RR$-space.
  \item
    A Stein manifold satisfies the {\em density property}
    if and only if it is an $RR$-space.
  \end{enumerate}
\end{conjecture}

For a first step towards $(i)$, regard Proposition~\ref{flex-tame}.

The last $RR$-axiom ($RR4$) is technically needed at several places.
For example, sometimes it is not enough to know that there is a very large
unstructured family of automorphisms, and one needs to have at least
one flow for a start. Is $(RR4)$ really necessary?

One step in removing $(RR4)$ might be:

\begin{conjecture}
  Let $X$ be a connected complex manifold on which its holomorphic
  automorphism group $\Aut(X)$ acts transitively.

  Then the connected component $\Aut^0(X)$
  of $\Aut(X)$ containing the identity map
  is non-trivial.
\end{conjecture}

\section{Linear Algebraic Groups}

\begin{theorem}\label{LAG-is-RR} 
  Let $G$ be a connected complex linear algebraic group
  which is {\em character-free}, i.e., there is no non-constant morphism
  of algebraic groups from $G$ to the multiplicative group $\C^*$.

  Then $G$ (as a complex manifold) is an $RR$-space.
\end{theorem}

\begin{proof}
  This follows from \cite{JW-LAG}, Theorem~1 and Theorem~3.
\end{proof}

\section{Preparations}

\subsection{Baire spaces}

A topological space is a {\em Baire space} if a countable intersection
of dense open sets is always dense. A set which is a
countable intersection of open
dense subsets is called {\em residual}.

For more information on Baire spaces see e.g. \cite{Munkres}, Chapter 8.

\begin{proposition}[=\cite{JW-LAG},Lemma~11]\label{baire}
  Let $X$ be a complex manifold.
  Then $\Aut(X)$ is a Baire space in the compact-open topology.
  \end{proposition}

The {\em compact open topology} may be generated
by the sets
  \begin{align*}
  &W(K,U)=\{\phi\in\Aut(X):\phi(K)\subset U\},\\
  &W'(K,U)=\{\phi^ {-1}\in\Aut(X):\phi(K)\subset U\}
  \quad\text{($K$ compact, $U$ open)}
  \end{align*}

  \begin{corollary}
    Let $X$ be a complex manifold.
    Let $\Aut(X)$ be equipped with the compact-open topology and let
    $\Aut^0(X)$ be the connected component containing the
    identity map.
  Then $\Aut^0(X)$ is a Baire space.
  \end{corollary}

  \begin{proof}
    As shown in \cite{JW-LAG} $\Aut(X)$ is  a Baire space because its
    topology may be defined by a complete metric.
    $\Aut^0(X)$ is closed in $\Aut(X)$ and closed subsets of complete
    metric spaces are again complete metric spaces.
  \end{proof}
  
\begin{proposition}[=\cite{JW-LAG}, Lemma~12]\label{aut-countable}
  Let $X$ be a complex manifold (as usual with countably topology).

  Then the topology of $\Aut(X)$ is likewise countable.
\end{proposition}

\begin{proposition}[=\cite{JW-LAG}, Corollary~8]\label{baire-1}
  Let $X$ be a complex manifold with a countable subset $D$ and an
  open dense subset $\Omega$. Assume that $\Aut(X)$ acts transitively
  on $X$.

  Then
  \[
  \{\phi\in\Aut(X):\phi(D)\subset \Omega\}
  \]
  is a {\em residual} subset of $\Aut(X)$.
  
  In particular, there exists an automorphism $\phi\in\Aut(X)$
  with $\phi(D)\subset\Omega$.
\end{proposition}

\begin{proposition}\label{baire-1a}
  Let $f$ be a non-constant holomorphic funciton on
  a complex manifold $X$ and let $D$ be a countable subset of $X$.
  Assume that $\Aut(X)$ acts $2$-transitively
  on $X$.

  Then 
  \[
  \{\phi\in\Aut(X):f\circ\phi|_D:D\to\C\text{ is injective}\}
  \]
  is a {\em residual} subset of $\Aut(X)$.
  
  In particular, there exists an automorphism $\phi\in\Aut(X)$
  such that $f$ is injective on $\phi(D)$.
\end{proposition}

\begin{proof}
  For every $x,y\in D$ with $x\ne y$ there exists an automorphism
  $\phi$ of $X$ such that $f(\phi(x))\ne f(\phi(y))$, because
  $f$ is non-constant and $\Aut(X)$ acts $2$-transitively.

  Hence for any $(x,y)\in D\times D$ with $x\ne y$
  the set
  \[
  \{\phi\in\Aut(X):f(\phi(x))\ne f(\phi(y))\}
  \]
  is open and dense.
  Since $D\times D$ is countable, it follows that
  \[
  \cap_{x\ne y\in D}\{\phi\in\Aut(X):f(\phi(x))\ne f(\phi(y))\}
  \]
  is residual.
\end{proof}

\begin{proposition}\label{baire-2}
  Let $X$ be a complex manifold. Assume that $\Aut(X)$ acts transitively
  on $X$. Let $\Omega\subset X$ be a non-empty Zariski-open subset
  (i.e., $X\setminus\Omega$ is a nowhere dense closed analytic subset).

  Then there exist finitely many elements $g_1,\ldots,g_r\in\Aut(X)$
  such that  
  \[
  X=\cup_{j=1}^rg_j(\Omega)
  \]
\end{proposition}

\begin{proof}
  For a given finite subset $S\subset\Aut(X)$, let
  \[
  X_S=\cup_{g\in S}g(\Omega)
  \]
  Let
  \[
  N=\min\{\dim(X\setminus X_S: S\subset\Aut(X), \text{$S$ finite}\}
  \]
  with $\dim\{\}\dfe-\infty$.
  \footnote{The dimension of a complex space is defined as the maximum
    of the dimensions of its irreducible components.}
  Assume that $N>-\infty$ and choose a finite subset $S_0\subset\Aut(X)$
  with
  \[
  N=\dim(X\setminus X_{S_0})
  \]
  Choose a countable subset $D\subset X$ such that $D$ intersects every
  connected component of $\Sigma\setminus \Sing(\Sigma)$ for
  $\Sigma=X\setminus X_{S_0}$.
  Due to Proposition~\ref{baire-1} there is an automorphism $\phi$ of $X$
  with $\phi(D)\subset X_{S_0}$.
  Set $S'=S_0\cup\{\phi\}$. By construction
  \[
  \dim(X\setminus X_{S'})<\dim(X\setminus X_{S_0})
  \]
  This contradicts minimality of $N$. Hence $N>-\infty$ may not occur.
  In other words, $N=-\infty$ and consequently
  \[
  X=\cup_{g\in S_0}X_{S_0}
  \]
\end{proof}

\subsection{Unbounded orbits}

\begin{lemma}\label{unbounded-orbits}
  Let $v$ be a complete vector field on a Stein space $X$.

  Then there are no relatively compact orbits of the associated
  flow except fixed points.
\end{lemma}

\begin{proof}
  Since $X$ is Stein, for every non-trivial orbit $\omega$ there exists
  a holomorphic function $f$ on $X$ which is not constant on $\omega$.
  Because every orbit is parametrized by
  a holomorphic map from $\C$, Liouvilles
  theorem implies that $f$ is unbounded on $\omega$.
  Therefore $\omega$ can not be
  relatively compact.
\end{proof}

\begin{lemma}
  Let $X$ be a Stein manifold on which its automorphism group $\Aut(X)$
  acts transitively, $v$ a complete vector field, $f$ a
  non-constant $v$-invariant
  function (i.e., $v(f)\equiv 0$),
  $\rho:X\to\R_0$ an exhaustion function and $D\subset X$
  an infinite discrete subset.

  Assume that  $\{x\in D:|f(x)|\le n\}$ is finite for every $n>0$, i.e.,
  $f|_D:D\to\C$ is proper.
  Assume in addition that  $f|_D$ is injective.

  Then for every function $\zeta:D\to\R^+$ there exists an automorphism
  $\psi$ of $X$ such that
  \[
  \rho(\psi(x))>\zeta(x)\ \forall x\in D
  \]
\end{lemma}

\begin{proof}
  There is an automorphism $\psi$ of $X$ such that $\psi(D)$
  does not intersect the zero set of $v$
  (Proposition~\ref{baire-1}).
  Thus there is no loss in generality
  in assuming that $v$ does not vanish at any point of $D$.
  Let $\phi_t$ denote the flow (one-parameter group of
  holomorphic automorphisms) associated to the vector field $v$:
  Then $D$ contains no fixed point of $\phi_t$ (because
  $v(x)\ne 0\ \forall x\in D$). Hence for every $x\in D$ the orbit of the
  flow $\phi_t$ through $x$ associated to $v$ is unbounded.
  (Lemma~\ref{unbounded-orbits}).
  It is therefore possible to choose a function
  $h_0:D\to\C$ such that
  \[
  \phi_{h_0(x)}(x)\not\in\{p\in X:\rho(p)\le \zeta(x)\}.
  \]
  Now we choose a holomorphic function $g:\C\to\C$ such that
  \[
  g(f(x))=h_0(x)\ \forall x\in D.
  \]
  Finally we set
  \[
  \psi:z\mapsto \phi_{g(f(z))}(z).
  \]
  This is an automorphism,
  because $f$ is invariant, which implies that $g\circ f$ is invariant
  and consequently $(g\circ f)v$ is a complete vector field.
  \end{proof}

\begin{corollary}\label{crit-tame}
  Let $X$ be a Stein manifold, $v$ a complete vector field, $f$ a
  non-constant invariant
  function  (i.e., $v(f)\equiv 0$), $D\subset X$
  an infinite discrete subset.

  Assume that $f|_D:D\to\C$ is proper and injective. 

  Then $D$ is weakly tame.
\end{corollary}

\begin{proposition}\label{g-w-tame}
  Let $X$ be a Stein manifold on which a connected complex Lie group $G$
  acts non-trivially.

  Let $D\subset X$ be an infinite discrete subset of $X$,
  $N\in\N$ and $F:X\to\C^N$ a holomorphic map such that
  \begin{enumerate}
  \item
    $D$ contains no fixed point for the $G$-action.
  \item
    $F$ is invariant under the $G$-action.
  \item
    $F$ restricts to a proper and injective map from $D$
    to $\C^N$.
  \end{enumerate}
  Then $D$ is weakly tame.
\end{proposition}

\begin{proof}
  Fix an exhaustion function $\rho$ on $X$ and a map
  $\zeta:D\to\R^+$.
  
  Recall that every non-trivial orbit is unbounded
  (Lemma~\ref{unbounded-orbits}).
  Hence $(i)$ implies that every point $x\in D$ is contained in
  an unbounded $G$-orbit.

  Thus there is a map $\xi:D\to G$ such that
  \[
  \rho(\xi(x)\cdot x)\ge\zeta(x)
  \]
  
  Further recall that there is a surjective holomorphic map $\alpha$ from some
  $\C^M$ to $G$.%
  \footnote{For example, this can be deduced from
    \cite{FF2017}, because every complex Lie group
    is an oka manifold.}
  We choose a map $\eta:D\to\C^M$ such that $\alpha\circ\eta=\xi$.

  Due to $(iii)$ the map $F$ maps $D$  bijectively onto a discrete
  subset $F(D)$ of $\C^N$.
  Hence there is a map $\beta:F(D)\to \C^M$ such that $\beta\circ F=\eta$.
  
  By standard complex analysis we may extend the map
  $\beta(D)\to\C^M$ to a holomorphic map
  $\gamma:\C^N\to\C^M$.

  We consider the map $\phi:X\to X$ defined as
  \[
  \phi(x)=\left(\alpha(\gamma(F(x)))\right)\cdot x
  \]
   
  Due to $(ii)$ we know that
  \[
  F(\phi(x))=F(x) \ \forall x\in X
  \]
  which implies that
  $\phi$ admits an inverse, namely
  \[
  x\mapsto \left(\alpha(\gamma(F(x)))\right)^{-1}\cdot x
  \]
  Thus $\phi\in\Aut(X)$.
  On the other, the construction implies that
  \[
  \forall x\in D:\rho(\phi(x))=\rho(\xi(x)\cdot x)\ge\zeta(x)
  \]

  This completes the proof that $D$ is weakly tame.
\end{proof}

  \begin{corollary}
  Let $X$ be a Stein manifold on which a connected complex Lie group $G$
  acts non-trivially. Assume that there exists a non-constant invariant
  (with respect to the $G$-action) holomorphic function $f$.

Then there exists a weakly tame discrete subset of $X$.
  \end{corollary}

  \begin{proof}
    We may replace $G$ by a one-dimensional subgroup and therefore
    assume that $G$ arises as the flow of a complete
    vector field
    $w$.
    Let $X^G$ denote the set of fixed points for the $G$-action
    on $X$, i.e., the zero locus of $w$.
    Denote $X^*=X\setminus X^G$. Consider the image
    $\Omega=f(X)$ which is a domain in $\C$.
    Note that $f(X^ *)$ is dense in $f(X)$.
    Hence we may find a sequence $q_n$ in $\Omega$ and
    a sequence $p_n\in X^ *$ such that:
    \begin{enumerate}
    \item
      $\{q_n:n\in\N\}$ is discrete in $\Omega$,
    \item
      $q_n\ne q_m$ for all $n\ne m$.
    \item
      $f(p_n)=q_n\ \forall n \in\N$.
    \end{enumerate}

    Let $D=\{p_n:n\in\N\}$. By construction the sequence
    $p_n$ has no convergent subsequence. Hence $D$ is discrete in $X$.

    We claim that $D$ is {\em weakly tame}.

    Let
    $\rho$ denote an exhaustion function on $X$ and
    let $\zeta:D\to\R^ +$ be a map.
    We have to show that there exists a holomorphic automorphism $\alpha$
    of $X$ with
    \[
    \rho(\alpha(p_n))>\zeta(p_n) \ \forall n\in\N
    \]
    
    Since $D\subset X^ *$, every $p_n\in D$ is contained in a
    one-dimensional $G$-orbit which is unbounded due to
    Lemma~\ref{unbounded-orbits}.
    We write the action as $\phi_t:X\to X$.
    For every $n\in\N$ we may choose a complex number
    $c_n\in\C$ such that
    \[
    \rho(\phi_{c_n}(p_n))>\zeta(p_n)
    \]
    Since $D'=\{q_n:n\in\N\}$ is a discrete subset of the domain $\Omega$,
    there is a holomorphic function $h\in{\mathcal O}(\Omega)$
    such that
    \[
    h(q_n)=h(f(p_n))=c_n\ \forall n \in \N
    \]
    Then $\tilde w=(h\circ f)w$ is a complete holomorphic vector
    field on $X$ and the associated flow $\tilde\phi$  satisfies
    \[
    \rho(\tilde\phi_1(p_n))=
    \rho(\phi_{c_n}(p_n))>\zeta(p_n)
    \]
  \end{proof}
    \section{First consequences of the axioms}
  
In this section $X$ always denotes a complex manifold which is
an $RR$-space.

\begin{proposition}\label{prop1}
  On an $RR$-space $X$ for a discrete subset $D$ the following conditions
  are equivalent:
  \begin{enumerate}
  \item
    $D$ is axiomatically tame.
  \item
    $D$ is weakly tame.
  \item
    $D$ is strongly tame.
  \end{enumerate}
\end{proposition}

\begin{proof}
  $(iii)\implies(ii)$ holds in general (see \cite{JW-TAME1}).
  
    $(i)\implies(iii)$ follows from axiom $(RR1)$ with $D=D'$.

  $(ii)\implies (i)$: Let $D$ be weakly tame. Let $\rho$ and $R_n$ be as in
  $(RR2)$. Recall that the choice of the exhaustion function is inessential
  (see \cite{JW-TAME1}). Hence we may choose the same exhaustion function for $(RR2)$
  as in the definition of the notion ``weakly tame''.
  Let $D=\{a_k:k\in\N\}$ and choose $\zeta:D\to\R^ +$
  with $\zeta(a_k)>R_k\ \forall k\in\N$.
  Since $D$ is weakly tame, there is an automorphism $\phi$ of $X$
  such that
  \[
  \rho(\phi(a_k))\ge\zeta(a_k)>R_k\ \forall k\in\N
  \]
  Then
  \[
  \#\{x\in\phi(D):\rho(x)\le R_n\}\le n\ \forall n\in\N
  \]
  Thus  $\phi(D)$ is axiomatically tame due to $(RR2)$.
  It follows that $D$ itself is axiomatically tame ($RR3$).
\end{proof}

In particular, on a given $RR$-space there is only one choice of
$\mathcal T$, and the notions of weak and strong tameness coincide.
Hence there is no ambiguity if from now on we simply speak about
``tame'' without specifying whether axiomatically, strongly or weakly
tame (as long as we are talking about infinite
discrete subsets in $RR$-spaces).

\begin{proposition}\label{prop2}
  Given any two tame discrete subsets $D$ and $D'$
  of an $RR$-space $X$ there is an
  automorphism $\phi$ of $X$ with $\phi(D)=D'$.
\end{proposition}

\begin{proof}
  Follows from $(RR1)$.
\end{proof}

\begin{corollary}
  The set $\mathcal T$ of tame discrete subsets is one
  $Aut(X)$-orbit in the set $\mathcal D$ of all
  infinite discrete subsets of $X$.
\end{corollary}

\begin{proposition}\label{unbounded-contains-tame}
  Every unbounded subset of an $RR$-space
  $X$ contains a tame discrete subset.
In particular, every infinite discrete set contains a tame discrete subset.
\end{proposition}

\begin{proof}
  Let $S$ be an unbounded subset. Then we choose a sequence $s_n\in S$ such that
  \[
  \forall n\in \N: \rho(s_{n+1})>\rho(s_n)>R_n
  \]
  Now axiom $(RR2)$ implies that $\{s_n:n\in\N\}$ is tame.
\end{proof}

\begin{proposition}\label{prop3}
  If $D$ is a tame discrete subset and $S$ is a finite subset, then
every infinite subset $D'$ of $S\cup D$ is tame.
\end{proposition}

\begin{proof}
  We may use induction on $\#S$.
  Hence it suffices to prove the case where $S$ consists of one point only,
  i.e.,
  $S=\{q\}$ for some $q\in X$.

  Choose $R_n$, $\rho$ as in axiom $(RR2)$ and let $D'$ be a discrete subset
  of $X$ satisfying
  \[
  \#\{x\in D':\rho(x)\le R_n\}=n-1\ \forall n\in\N
  \]
  $(RR2)$ implies that $D'$ is tame. Since $D$ and $D'$ are both tame,
  there is a holomorphic automorphism $\phi$ of $X$ with $\phi(D)=D'$
  (Proposition~\ref{unbounded-contains-tame}).
  Let $q'=\phi(q)$ and $D''=D'\cup\{q'\}$.
  Now
  \[
  \#\{x\in D'':\rho(x)\le R_n\}\le n\ \forall n\in\N
  \]
  and $D''$ is tame due to $(RR2)$.
  Therefore $D\cup\{q\}=\phi^{-1}(D'')$ is tame, too.
  \end{proof}

\begin{corollary}\label{sym-diff}
  Let $D$, $D'$ be infinite discrete subsets such that
  their ``symmetric difference'', i.e.,
$(D\setminus D')\cup(D'\setminus D)$
is finite. Then $D$ is tame if and only if $D'$ is tame.
\end{corollary}

\begin{proof}
  $D'$ is an subset of $D\cup(D'\setminus D)$ and $D'\setminus D$
  is finite. Hence Proposition~\ref{prop3} implies that $D'$ is tame
  if $D$ is tame. The opposite direction follows similarily.
\end{proof}

\begin{corollary}\label{cor-sub-tame}
  Every finite subset is contained in a (infinite)  tame discrete
  subset.
\end{corollary}

\begin{proof}
  Let $S$ be a finite subset. Due to $(RR4)$ (or $(RR2)$) there exists a tame
  discrete subset $D$.
  Then $S\cup D$ is tame due to Proposition~\ref{prop3}.
\end{proof}

\begin{proposition}\label{inf-trans}
  ($\infty$-transitivity.)

  Let $X$ be an $RR$-space and let $D\subset X$ be empty or a tame discrete
  subset.
  
  Let $S$, $S'$ be finite subsets of $X\setminus D$
  and let $f:S\to S'$ be an injective
  map. Then $f$ extends to a holomorphic automorphism $\phi\in Aut(X)$
  with $\phi|_D=id_D$.
\end{proposition}

{\em Remark.} For finite sets of cardinality $n$ this property is
called $n$-transitivity of the $\Aut(X,D)$-action on $X\setminus D$.
($\Aut(X,D)$ denotes the set of holomorphic automorphisms $\phi$ of $X$
with $\phi(D)=D$.)

\begin{proof}
  If $D$ is empty, we replace $D$ by an arbitrary infinite tame discrete
  subset of $X$ with $D\cap(S\cup S')=\{\}$.
  
  Both $D_1=D\cup S$ and $D_2=D\cup S'$ are tame (Proposition~\ref{prop3}).
  Let $f_1:D_1\to X$ be defined as $f_1|_D=id_D$ and $f_1|_S=f$.
  Now $f_1:D_1\to D_2$ is an injective map
  which extends to an automorphism $\phi$ due to axiom $(RR1)$.
  By construction we have $\phi|_D=id_D$ and $\phi|_S=f$.
\end{proof}

\begin{proposition}\label{no-equ-rel}
  Let $\sim$ be an equivalence relation on $X$ which is invariant
  under all automorphisms.

  Then $\sim$ is trivial, i.e.,
  either 
  \[
  x\sim y\ \ \forall x,y\in X
  \]
  or
  \[
  x\sim y \iff x=y.
  \]
\end{proposition}

\begin{proof}
  Assume that there exist $x,y\in X$ with $x\sim y $, but $x\ne y$.
  Let $p,q\in X$ be arbitrary points with $p\ne q$.
  Now $Aut(X)$ acts ``$2$-transitively'' by Proposition~\ref{inf-trans}.
  Thus there is an automorphism $\phi$ with $\phi(x)=p$ and $\phi(y)=q$.
  Since $\sim$ is supposed to be $Aut(X)$-invariant, we may conclude that
  $p\sim q$ for all $p,q\in X$.
\end{proof}

\begin{proposition}\label{no-metric}
  On a $RR$-space $X$ there is no non-trivial pseudometric
  invariant under all automorphisms.

  In particular both the Kobayashi and the Caratheodory
  pseudo distance vanish. As a consequence, every bounded
  holomorphic function on $X$ is constant.
\end{proposition}

\begin{proof}
  A non-trivial pseudometric $\delta(\ ,\ )$ defines a non-trivial
  equivalence relation via
  \[
  x\sim y \quad\iff\quad \delta(x,y)=0
  \]
  Hence the assertion follows from Proposition~\ref{no-equ-rel}.
\end{proof}

\begin{proposition}\label{stein}
Every $RR$-space is Stein.
\end{proposition}

\begin{proof}
  Due to axiom $(RR4)$ there is a non-constant holomorphic function
  $f$ on $X$.

Since the equivalence relation defined by holomorphic functions is
trivial (Proposition~\ref{no-equ-rel}),
$X$ must be holomorphically separable, i.e.,
for every $x,x'\in X$ with $x\ne x'$ there is a holomorphic function
$F:X\to\C$ with $F(x)\ne F(x')$.

Proposition~\ref{no-metric} implies that our
non-constant holomorphic function $f$ is unbounded.

Choose a sequence $p_n\in X$ such that
\[
\exists N:\forall n\ge N:\ n<|f(p_n)|<n+1
\]
(this is possible, because the image $f(X)$ is unbounded
and connected) and define
$E=\{p_n:n\in\N\}$. Let $D$ be any
discrete set. 
Then $D$ and $E$ contain tame discrete subsets $D_1$ resp.~$E_1$
(proposition~\ref{unbounded-contains-tame}).
Let $\phi$ be an holomorphic automorphism of $X$ with $\phi(D_1)=E_1$
(by axiom $(RR1)$).
Then $f\circ\phi$ is unbounded on $D$.
Thus we proved: For any infinite discrete subset $D$ of $X$ there
exists a holomorphic funciton on $X$ which is unbounded on $D$.
Therefore $X$ is
holomorphically convex.

Thus $X$ is both holomorphically separable and holomorphically convex
and therefore Stein (\cite{GR},IV,\S12, Fundamental Theorem).
\end{proof}

\begin{proposition}\label{not-loc-cpt}
  The automorphism group $Aut(X)$ of a $RR$-space
  $X$ is not locally compact.
\end{proposition}

\begin{proof}
  Assume that $\Aut(X)$ is locally compact. Because the
  topology of $\Aut(X)$ is countable
  (Proposition~\ref{aut-countable}), this implies that there
  are countably many compact subsets $K_n$ of $\Aut(X)$
  with $\cup_{n\in\N}K_n=\Aut(X)$.
  
  Let $D=\{p_n:n\in\N\}$ be a tame discrete subset of $X$.
  We choose an injective map $f_0:D\to D$ such that
  \[
  \forall n\in\N: f_0(p_n)\not\in K_n\cdot p_n=\{\psi(p_n):\psi\in K_n\}
  \]
  (This is possible, because $D$ is unbounded and $K_n\cdot p_n$
  is compact for every $n$.)
  Due to axiom $(RR1)$ there must exist a holomorphic automorphism
  $\phi\in\Aut(X)$ with $\phi(p_n)=f_0(p_n)\ \forall n\in\N$.
  Since $\Aut(X)=\cup_n K_n$, there is an index $n\in\N$ such that
  $\phi\in K_n$. This yields a contradiction:
  On one side, $\phi(p_n)=f(p_n)$, on the other side
  $f(p_n)\not\in K_n\cdot p_n$.
\end{proof}

\begin{proposition}\label{non-tame}
  Every $RR$-space admits a {\em non-tame} infinite discrete
  subset.
\end{proposition}

\begin{proof}   Since every $RR$-space is Stein (Proposition~\ref{stein}),
   theorem~2 of
   \cite{JW-LDS} implies that there is a discrete subset $D$
   such that $\phi|_D=Id_D$ for every holomorphic automorphism $\phi$
   of $X$ with $\phi(D)=D$.
   But condition $(RR1)$ implies in particular that every bijective self-map
   of a tame discrete subset extends to a holomorphic automorphism of $X$.
   Hence $D$ can 
   not be tame.
\end{proof}

\begin{proposition}
  Let $X$ be an $RR$-space, $Z$ a closed subset, $v$ a complete
  vector field with flow $\phi$
  and $f$ an $\phi$-invariant holomorphic
  function (i.e.~$v(f)=0$).

  Assume that $f|_Z:Z\to\C$ is an unbounded proper continuous map.

  Then an infinite discrete subset $D\subset X$ is tame if and only
  if there exists an automorphism $\phi\in\Aut(X)$
  with
  $\phi(D)\subset Z$.
\end{proposition}

\begin{proof}
  $Z$ is not compact, because $f|_Z$ is unbounded.
  Therefore $Z$ contains a tame discrete subset $D_0$
  (proposition~\ref{unbounded-contains-tame}).
  For every tame discrete subset $D$ of $X$ there is an
  automorphism $\phi$  of $X$ with $\phi(D)=D_0$,
  implying $\phi(D)\subset Z$.

  Conversely, let $D\subset Z$ be an infinite discrete subset.
  Consider the set
  \[
  T=\{\phi\in\Aut(X):\forall x\in D:\phi(x)\in Z\}.
  \]
  This is a closed subset of $\Aut(X)$ and therefore a Baire space.
  Let $p\in D$ and $c\in Z$. Then $D'=(D\cup\{c\})\setminus p$
  is tame (Corollary~\ref{sym-diff}) which implies that every
  bijection from $D$ to $D'$ extends to an automorphism of $X$
  ($RR1$).
  Given $p,q\in D$ with $p\ne q$ we may choose $c\in Z$ with
  $f(c)\ne f(q)$ (because $f|_Z$ is not constant) and therefore
  find an automorphism $\phi\in T$ with
  \[
  f(\phi(p))=f(c)\ne f(q)=f(\phi(q))
  \]
  It follows that
  \[
  T_{p,q}=\{\phi\in T: f(\phi(p))\ne f(\phi(q))\}
  \]
  is an open dense subset of $T$ and that
  \[
  \cap_{p,q\in D; p\ne q} T_{p,q}
  \]
  is a residual subset of $T$ and therfore not empty.
  This implies that there is an auotmorphism $\phi$ of $X$
  such that
  \begin{enumerate}
  \item
    $\phi(D)\subset Z$.,
  \item
    $f\circ\phi|_D:D\to\C$ is injective.
  \end{enumerate}
  
  Now tameness of $D$ follows from Corollary~\ref{crit-tame}.
\end{proof}

\subsection{Axiom $(RR4)$ and Generalized shears/replicas}

A very important class of examples for the axiom $(RR4)$ is
related to {``\em shears''}:

Let $X=\C^n$ with $v=\frac{\partial}{\partial z_{1}}$
and $f$ being any function which depends only on $z_2,\ldots,z_n$
and not on $z_1$.

Then $(v,f)$ satisfy the conditions of axiom $(RR4)$.

The flow $\phi_t$ associated to $v$ is given as
\[
\phi_t(z_1,\ldots,z_n)\mapsto (z_1+t,z_2,\ldots,z_n)
\]

Now $fV$ is again a complete vetcor field, with associated
flow
\[
\psi_t:(z_1,\ldots,z_n)\mapsto (z_1+tf(z_2,\ldots,z_n),z_2,\ldots,z_n)
\]
Such
an automorphism $\psi_t$ of $\C^n$ is often called a {``\em shear''},
the vector field $fv$ is often called (in particular in affine
algebraic geometry considering locally nilpotent derivations)
a {``\em replica''} of $v$.

This construction generalizes as follows:

\begin{proposition}
  Let $X$ be a complex manifold, $v$ a complete vector field
  with flow $\phi_t$ and let $f$ be a holomorphic function on $X$
  with $v(f)=0$.

  Let $h:\C\to\C$ be an arbitrary entire function.

  Then
  \begin{enumerate}
  \item
    $(h\circ f)v$ is a complete vector field.
  \item
    The flow associated to $(h\circ f)v$ may be described as
    \[
    \xi_t:x\mapsto \phi_{h(f(x))t}(x)
    \]
  \end{enumerate}
\end{proposition}
  
This process, where we replace a complete vector field by itself multiplied
with some invariant function, may be visualized as keeping the orbit
structure but modifying the speed with which the flow moves through
the orbits.

\begin{proposition}
  Let $X$ be a Stein complex manifold with a complete vector field $v$ and an
  invariant (i.e. $v(f)=0$) holomorphic function $f$.

  Let $D$ be a discrete subset of $X$ such that the restricted map
  $f|_D:D\to\C$ is proper.
  Assume that $v_p\ne 0\ \forall p\in D$.

  Assume in addition that
  \begin{enumerate}
  \item
    $f|_D$ is injective, or
  \item
    all the orbits of the flow associated to $V$ are closed.
  \end{enumerate}

  Then $D$ is weakly tame, i.e., $D$ satisfies $(RR2)$.
\end{proposition}

\begin{remark}
  In the setup of the proposition,
  assume that $X$ is an affine algebraic variety,
  that the flow associated to $v$
  is an algebraic action of the additive group $(\C,+)$.

  Then all the orbits are necessarily closed, thus $(ii)$ is satisfied.
\end{remark}

\begin{proof}
  For every value $q\in f(D)$, the preimage $\{x\in D:f(x)=q\}$ is finite.
  Recall that the flow $\phi_t$ associated to $v$ has only unbounded
  orbits except fixed points
  (Lemma~\ref{unbounded-orbits}). However, since we assumed
  $v_p\ne 0\ \forall p\in D$, no element of $D$ is a fixed point.
  Hence every orbit through a point of $D$ is unbounded.
  
  Enumerate $f(D)$ as $f(D)=\{a_k:k\in\N\}$ and define
  \[
  \sigma(k)=\sum_{j\le k}\#f_D^{-1}(a_j)
  \]
  Let $\rho$, $R_n$ be as provided by axiom $(RR2)$. Since $\rho$
  is an exhaustion function, $\{x\in X:\rho(x)\le R_n\}$ is
  compact for every $n$.

  $(i)$ Assume that $f|_D$ is injective.
  Enumerate $D$ as $\{c_k:k\in\N\}$ such that $f(c_k)=a_k$.
  Then we choose
  $\alpha_k\in\C$ such that
  \[
  \phi_{\alpha_k}(c_k)\not\in\{x:\rho(x)\le R_k\}
  \]
  (which is possible, because the orbit of the flow through $c_k$ is
  unbounded.)

  $(ii)$ Here we assume that the orbits of the flow are closed.
  In this case,
  \[
  \{x\in\omega:\rho(x)\le R_k\}
  \]
  is compact for every orbit $\omega$ and every $k\in\N$.
  Hence for each $k$ we may choose $\alpha_k\in\C$ such that
  \[
  \phi_{\alpha_k}(c)\not\in\{x:\rho(x)\le R_{\sigma(k)}\}
  \ \forall c\in D, f(c)=a_k.
  \]

  In both cases our construction implies
  \[
  \#\left\{x\in D: \exists k: f(x)=a_k,
  \rho\left(\phi_{\alpha_k}(x)\right)\le R_n
  \right\}\le n \quad  \forall n
  \]
  We choose a holomorphic function $h:\C \to\C$ with
  $h(a_k)=\alpha_k\ \forall k\in\N$ and define $w$ as the
  vector field $(h\circ f)v$. Then $w$ is complete and its time $1$
  map $\psi$ satisfies
  \[
  \#\left\{x\in\psi(D):\rho(x)\le R_n\right\}\le N.
  \]
  Hence $D$ is weakly tame, i.e., satisfies $(RR2)$.
\end{proof}

\begin{proposition}\label{embed}
  Let $G$ be a group acting 2-transitively on a complex manifold $X$
  and let $\pi:X\to Y$ be a non-constant holomorphic map between complex
  manifolds.

  Then there exist finitely many elements $g_1,\ldots,g_k\in G$
  such that
  the map $\Phi:X\to Y^k$ defined as
  \[
  \Phi:x\mapsto (\pi(g_1(x)),\ldots,\pi(g_k(x)))
  \]
  is injective and immersive.
\end{proposition}

\begin{proof}
  For every finite subset $S=\{g_1\ldots,g_r\}\subset G$ let
  $
  \Phi_S:X\mapsto Y^S$
  be defined as
  \[
  \Phi_S:x\mapsto (\pi(g_1(x)),\ldots,\pi(g_k(x)))
  \]
  and let $d_S$ denote the minimal dimension of irreducible components
  of fibers of $\Phi_S$.

  We claim that there is a finite set $S$ with $d_S=0$. Assume the
  contrary, let $d_0$ be the minimal possible value of $d_S$ and let
  $S_0\subset G$ be a finite subset with $d_0=d_{S_0}$.
  Let $\Sigma$ be a $d_0$-dimensional irreducible component of
  a fiber $X_p=\Phi_{S_0}^{-1}(p)$ of $\Phi_{S_0}:X\to Y^{S_0}$,
  and choose two distinct  points
  $q,\tilde q\in \Sigma\setminus\Sing(X_p)$.%
  \footnote{Here $\Sing(X_p)$ denotes the singular locus of the
    {\em reduced fiber}. Hence $\Sing(X_p)\cap\Sigma$ is nowhere dense
    in $\Sigma$.}
  Since $G$ acts $2$-transitively on $X$, we may choose an element
  $g\in G$ such that $g(q)=q$ and $g(\tilde q)\not\in X_p$.
  Let $S'=S_0\cup\{g\}$. By construction, the $\Phi_{S'}$-fiber through
  $q$ has dimension less than $d_0$. Therefore $d_{S'}<d_{S_0}=d_0$, which
  contradicts the minimality of $d_0$.

  Therefore $d_0=0$ and consequently there is a finite set $S_0$ and a point
  $q\in X$ such that $q$ is an isolated point in the $\Phi_{S_0}$-fiber
  through $q$.
  By semicontinuity of fiber-dimension for holomorphic maps it follows
  that there is an open subset $W\subset X$ such that every $w\in W$ is
  an isolated point in the $\Phi_{S_0}$-fiber through $w$.
  Due to Sard's theorem we may conclude that there is a dense subset of
  $W$ where $\Phi_{S_0}$ is immersive.
  The set $\Omega$ of all points in $X$ where $\Phi_{S_0}$ is immersive
  is a Zariski open subset.
  Hence using Proposition~\ref{baire-2} there are finitely many
  elements $g_1,\ldots,g_s\in G$ such that
  \[
  X=\cup_{j=1}^sg_j(\Omega)
  \]
  Let
  \[
  S_1=\cup_{j=1}^s\cup_{\alpha\in S_0}g_j\cdot\alpha
  \]
  Then $\Phi_{S_1}$ is everywhere immersive.

  In a similar spirit, we deal with injectivity.
  
  Let
  \[
  P=\{(x,y)\in X\times X,x\ne y\}/\sim\quad\text{$(x,y)\sim(y,x)$}
  \]
  Since the $G$-action on $X$ is $2$-transitive, $G$ acts transitively on
  $P$.
  For any finite set $S\subset G$
  we consider the {\em incidence variety}
  \[
  Z_S=\{[(x,y)]\in P:\Phi_S(x)=\Phi_S(y)\}
  \]
  Let $D_S=\dim(Z_S)$, i.e., $D_S$ is the maximal dimensional of
  irreducible components of $Z_S$.
  Fix $S$ such that $D_S$ is minimal and choose a countable subset
  $\Lambda\subset Z_S$ which intersects the smooth locus of every
  irreducible component of $Z_S$.
  Proposition~\ref{baire-2} implies the existence of an automorphism
  $\alpha\in G$ such that
  \[
  \alpha(\Lambda)\subset P\setminus Z_S
  \]
  It follows that for $S'=S\cup \alpha S$ we have
  \[
  Z_{S'}\subset Z_S\text{ and } Z_{S'}\cap\Lambda=\{\}
  \]
  By construction of $\Lambda$ it follows that $\dim(Z_{S'})<\dim(Z_S)$.
  This yields a contradiction to minimality unless $D_S=-\infty$, i.e.,
  unless $\Phi_S$ is injective.

  Finally, if $\Phi_S$ is injective and $\Phi_{S'}$ is immersive,
  then $\Phi_{S\cup S'}$ is both injective and immersive.
\end{proof}

\begin{lemma}\label{L2.16}
  Let $X$ be a complex manifold
  with a non-zero complete vector field $v$.
  Let $G\subset\Aut(X)$  be a subgroup.
  
  Assume that $G$ acts $2$-transitively on $X$.

  Then:
   {\em Every tangent space
   $T_pX$ is generated as a vector space by
   the tangent vectors $(\phi^*v)_p$ with $\phi\in G$. }
\end{lemma}

\begin{proof}
  Let $\cH$ be the subsheaf of the tangent sheaf $\cTX$ which
  is the ${\mathcal O}_X$-module sheaf
  generated by all vector fields of the form $\phi^*v$ with $\phi\in G$.

  Since $v$  is a non-zero vector field, $\cH$
  is not the zero sheaf. On the other hand, $\cH$ is invariant under all
  automorphisms $\phi\in G$ and $G$ acts transitively on $X$.
  Thus $\cH$ is locally free.

  Let $\chH$ denote the smallest involutive subsheaf generated by
  $\cH$, i.e., $\chH$ is generated by the elements of the Lie algebra
  generated by (local) sections of $\cH$. Then $\chH$ defines a foliation.

  We claim that $\chH=\cTX$.
  Let $H$ denote the subgroup of $\Aut(X)$ generated by all the flows of
  $\phi^*v$ for all $\phi\in G$. Note that for a fixed $\phi$ the flow
  of $\phi^ *v$ is conjugate by $\phi$ to the flow of $v$.
  As a consequence, $H$ is normalized $G$.
  Thus, if $\hat G$ denotes the subgroup of
  $\Aut(X)$ generated by $H$ and $G$, tghen $H$ is a normal subgroup
  of $G$.
  
  Choose a point $p$ at which $v$ is not zero. Then the flow of $v$ is
  not contained in the isotropy group
  \[
  J_p=\{\phi\in \hat G:\phi(p)=p\}
  \]
  Hence $H\not\subset J_p$ and therefore $J_p\not\subset H\cdot J_p$.
  Recall that $G$ is assumed to act $2$-transitively on $X$.
  It follows $\hat G$ acts $2$-transitively on $X=\simeq \hat G/J_p$
  This implies that there does not exist any subgroup $I$ with
  \[
  J_p\subsetneq I \subsetneq \hat G
  \]
  Since $J_p\not\subset H\cdot J_p$, we may now deduce $H\cdot J_p=\hat G$.
  In particular, $H$ acts transitively on $X$. But the construction of
  $H$ implies that the $H$-orbits in $X$ are contained in the leaves of
  the  foliation defined by $\chH$.
  It follows that $\chH=\cTX$.

   Next we claim that  $\cH$ itself is already involutive,
  i.e., $\cH=\chH$.

  To verify this, assume the contrary.
  By construction of $\cH$ this means that there are
  finite index sets $I,J$, elements $\alpha_i,\beta_j\in G$,
  a point $p\in X$
  and holomorphic functions $f_i$, $g_j$ near $p$
  such that
  \[
    [u,w]_p\not\in\cH_p
    \]
    for  $u=\sum_i f_i\alpha_i^*v$, $w=\sum_j g_j\beta_j^*v$.
  Hence there exist $i_0\in I, j_0\in J$ with
  \[
    [f_{i_0}\alpha_{i_0}^*v,g_{j_0}\beta_{j_0}^*v]_p\not\in\cH_p
  \]
  Observe that
  \[
    [f_{i_0}\alpha_{i_0}^*v,g_{j_0}\beta_{j_0}^*v]
    =
    \underbrace{f_{i_0}\alpha_{i_0}^*v(g_{j_0})\beta_{j_0}^*v-g_{j_0}\beta_{j_0}^*v(f_{i_0})\alpha_{i_0}^*v}%
    _{\in \cH}
    +f_{i_0}g_{j_0}[\alpha_{i_0}^*v,\beta_{j_0}^*v].
  \]
    
  Thus
  \[
    [\tilde u,\tilde w]\not\in\cH_p
  \]
  for $\tilde u =\alpha_{i_0}^*v$ and $\tilde w=\beta_{j_0}^*v$.
  \renewcommand{\vv}{\tilde u}
  \renewcommand{\ww}{\tilde w}
  Let $\psi_t$ denote the flow of $\vv$ and recall that $[\vv,\ww]$ equals
  the Lie derivative of $\ww$ with respect to $\vv$.
  Thus
  \[
    [\vv,\ww]_p=\lim_{t\to 0}\frac{\left(\psi_t^ *\ww\right)-\ww}{t}.
   \]
   We assumed $\phi_t\in G$. Now
   $\psi_t=\alpha_{i_0}\circ\phi_t\circ(\alpha_{i_0})^{-1}$,
   because $\vv=\alpha_{i_0}^*v$.
   Hence $\psi_t\in G$ and consequently both $\phi_t^*\ww$ and
   $\ww$ are sections
  in $\cH$.
  This yields a contradiction to $[\vv,\ww]\not\in\cH_p$.
 
  Thus $\cH$ defines a foliation.
  It follows that $\cH=\chH=\cTX$.
  \end{proof}

\begin{corollary}\label{H-trans}
  Let $X$ be a complex manifold
  with a non-zero complete vector field $v$.
  Let $G\subset\Aut(X)$  be a subgroup.
  
  Assume that $G$ acts $2$-transitively on $X$.

  Let $H$ denote a subgroup of $\Aut(X)$
  containing the flow of $\phi^*v$ for all $\phi\in G$.
  
  Then $H$ acts transitively on $X$.
\end{corollary}

\begin{proof}
  The Lie algebra of $H$-fundamental vector fields on $X$ contains
  $\phi^*v$ for all $\phi\in G$. Hence the $H$-fundamental
  vector fields span $TX$ at every point of $X$ by Lemma~\ref{L2.16},
  (ii). It follows that $H$ acts transitively on $X$.
\end{proof}

\begin{proposition}\label{H-n-trans}
  Let $X$ be a complex manifold
  with a non-zero complete vector field $v$.
  Let $I=\{\phi_t:t\in\C\}$
  denote the one-parameter-group of automorphisms of $X$
  arising from the flow of $v$.
  Let $G$ be a subgroup of $\Aut(X)$
  which acts $2n$-transitively on $X$ ($n\in\N$).
    
  Let $H$ denote a normal subgroup of $G$
  containing $I$.
  
  Then $H$ acts $n$-transitively on $X$.

\end{proposition}

\begin{corollary}\label{H-infi-trans}
  In the setup of Proposition~\ref{H-n-trans}
  assume in addition that  $\Aut(X)$ acts $\infty$-transitive on $X$.
  
  Then $H$ acts $\infty$-transitively on $X$.
\end{corollary}

\begin{proof}[Proof of \ref{H-n-trans}]
  Let $Y$ denote the open subset of $X^n$ containing those
  $(x_1,\ldots,x_n)$ which are all distinct ($x_i\ne x_j$ if $i\ne j$).
  Note that $v$ induces a complete vector field $\tilde v$ on $Y$
  with associated flow
  \[
  \tilde\phi_t:(x_1,\ldots,x_n)\mapsto
  \left(\phi_t(x_1),\ldots,\phi_t(x_n)\right)
  \]
  By assumption $G$ acts transitively on $Y$.
  Fix an element $x=(x_1,\ldots,x_n)\in Y$ and let $J$ denote the
  isotropy of $G$ at $x$.
  Since the $G$-action on $X$ is $2n$-transitively,
  $J$ acts transitively on the open set
  \[
  \Omega=\{(y_1,\ldots,y_n)\in Y:y_j\ne x_j\forall j\}
  \]
  Furthermore, for every $g\in G$ the flow associated to
  $g^*v$ is conjugate to the flow of $v$ via $g$. Since $H$ is
  assumed to be normal, it follows that $H$ contains the
  flow of $g^*v$ for every $g\in G$.
  
  We may thus apply Corollary~\ref{H-trans} and conclude
  that $H$ acts transitively on $Y$.
  It follows that the $H$-action on $X$ is $n$-transitive.
\end{proof}

\begin{corollary}\label{0-inf-trans}
  Let $X$ be an $RR$-space.

  Then the connected component $\Aut^0(X)$ of $\Aut(X)$ which contains
  the identity self-map acts $\infty$-transitively on $X$.
\end{corollary}

\begin{proof}
  Every $RR$-space supports a complete vector field (Axiom $(RR4)$).
  The assertion follows from
  Proposition~\ref{inf-trans} and Corollary{H-infi-trans}
  with $H=\Aut^ 0(X)$.
\end{proof}

\begin{proposition}\label{t-span-vf}
  Let $X$ be an $RR$-space. Let $D\subset X$
  be finite (possibly empty) or a tame infinite discrete subset.

  Let $X^*=X\setminus D$ and
  $G=\{\alpha\in\Aut(X): \alpha(x)=x\ \forall x\in D\}$.

  Let $G^0$ denote the connected component of $G$ containing the
  neutral element.
  
  \begin{enumerate}
  \item
    $G^0$ acts infinitely transitive on $X^*$.
  \item
    The tangent bundle of $X^*$ is everywhere spanned by
    complete vector fields $v$ which vanish on $D$.
  \end{enumerate}
\end{proposition}

\begin{proof}
  Due $(RR4)$ the $RR$-space $X$ admits a non-zero complete vector field
  $v$ and a non-constant  holomorphic function $f$ with $v(f)\equiv 0$.

  We note that $f$ is unbounded, since every bounded holomorphic function
  on $X$ is constant due to Proposition~\ref{no-equ-rel}.

  In combination with $(RR2)$ this allows use to find a tame discrete
  subset $D'$ such that $f(D')$ is discrete in $\C$.
  We choose a non-constant holomorphic function $h$ on $\C$ which vanishes
  on $f(D')$ and define $F=h\circ f$.  

  By Proposition~\ref{prop2}  there is a holomorphic automorphism $\phi$ of $X$
  with $\phi(D)=D'$.

  Define $\tilde v=Fv$ and $w=\phi^*\tilde v$.
  Now $\tilde v$ is a complete vector field on $X$, because $v$ is a complete
  vector field and $v(f)\equiv 0$ implies
  $v(F)=v(h\circ f)\equiv 0$.
  
  By construction $F(\phi(x))=0\ \forall x\in D$. It follows that
  $w_p=0$ for all $p\in D$. Thus $w$ is complete on $X^*=X\setminus D$.
  Evidently the flow of $w$ is contained in $G^ 0$.
  Note that $G^0$ is a normal subgroup of $G$.

  Proposition~\ref{inf-trans} implies that $G$ acts $\infty$-transitively on
  the manifold $X^*$.
  Thus claim $(ii)$ follows from
  Lemma~\ref{L2.16}, $(ii)$.

  Claim $(i)$ follows from
  Corollary~\ref{H-infi-trans}.
\end{proof}

\begin{proposition}\label{span-finite}
  Let $X$ be a complex manifold. Assume that the tangent bundle is
  everywhere spanned by complete vector fields.

  Then there are {\em finitely} many complete vector fields which
  everywhere span the tangent bundle.
\end{proposition}

\begin{proof}
  Since complete vector fields generate the tangent bundle everywhere,
  $\Aut(X)$ acts transitively on $X$.
  Choose $p\in X$ and let $v_1,\ldots,v_n$ be complete vector fields
  generating the tangent space $T_pX$. Denote $A_1=\{v_1,\ldots,v_n\}$.
  Let $Z_1$ denote the subset of $X$ where $v_1,\ldots,v_n$ fail to
  generate the tangent space. This is a closed analytic subset with
  countably many irreducible components $S_{1,j}$. We choose a point
  $p_{1,j}$ in the smooth part of $S_{1,j}$ for every $j$.
  For each $j$, let $\Omega_j$ denote the set of all $g\in\Aut(X)$ such that
  $g_*v_j$ span $T_{p_{1,j}}X$. Then $\Omega_j$ is open and dense in $\Aut(X)$.
  Now $\Aut(X)$ is a Baire space (Proposition~\ref{baire}).
  Hence $\cap_j \Omega_j$ is not empty and we may choose
  $g_1\in\cap_j \Omega_j$.
  Let $A_2=A_1\cup\{(g_1)_*v_j\}$ and let $Z_2$ denote the set of all points
  of $X$ where the elements of $A_2$ fail to generate the tangent space.
  By construction
  $Z_2$ is a nowhere dense analytic subset of $Z_1$. In particular,
  $\dim(Z_2)<\dim(Z_1)$.
  Proceeding recursively, we find
  a sequence of closed analytic subsets with
  $\dim(Z_{k+1}<\dim(Z_k)\ \forall k$.
  Evidently this sequence must terminate.
  This yields a finite family of complete
  vector fields which generate the tangent space everywhere.
\end{proof}

\begin{remark}
  In the above proposition it is essential that we consider
  consider the set of {\em all}
  complete vector fields.

  For example, let $f_n:\C\to\C$ be an entire function
  with $\{z\in\C:f_n(z)=0\}=\{k\in\N:k\ge n\}$.

  Then $T\C^2$ is everywhere generated by the vector fields
  \[
  \frac{\partial}{\partial z_1},\quad f_n(z_1)\frac{\partial}{\partial z_2}
    \quad(n\in\N)
    \]
    but it is not generated by any finite subset
    of these complete vector fields.
\end{remark}

\begin{proposition}\label{oka}
  Let $X$ be an $RR$-space and $D$ a discrete subset of $X$ which is
  empty, finite or tame.

  Then $X^*=X\setminus D$ is an Oka manifold.
\end{proposition}

\begin{proof}
  Due to Proposition~\ref{t-span-vf} and Proposition~\ref{span-finite}
  the tangent bundle of $X^*$ is spanned by finitely many complete vector
  fields $v_1,\ldots,v_N$ on $X^*$.
  Let $\phi_{j,t}$ denote the corresponding flows.
  Let $E=X^*\times\C^N$ be the trivial vector bundle and
  define
  \[
  s:(x;t_1,\ldots,t_N)=\phi_{1,t_1}\circ\ldots
  \phi_{N,t_N}(x)
  \]
  This is a spray and therefore $X^*$ is elliptic
  and in particular Oka.
\end{proof}

\begin{corollary}\label{cn-surj}
  For every $RR$-space there exists a surjective holomorphic map
  from some $\C^N$ onto $X$.
\end{corollary}

\begin{proof}
  Follows from Proposition~\ref{oka} and \cite{FF2017}.
\end{proof}

\begin{corollary}\label{no-bd-psh}
  A bounded plurisubharmonic function $\rho$ on a $RR$-space $X$
  is constant.
\end{corollary}

\begin{proof}
  Let $F:\C^N\to X$ be a surjective holomorphic map.
  Then $\rho\circ F$ is a bounded plurisubharmonic
  function on $\C^N$ and therefore constant.
  Thus $\rho$ itself must be constant.
\end{proof}

\begin{proposition}\label{avoid}
  Let $X$ be an $RR$-space and let $D$ be tame discrete subset.
  Let $Y$ be a Stein complex manifold with $\dim(X)=\dim(Y)$.
  Then there exists a holomorphic map $F:Y\to X\setminus D$
  which has full rank at some point.
\end{proposition}

\begin{proof}
  Let $n=\dim(X)$.
  $X\setminus D$ is an Oka manifold (Proposition~\ref{oka}).
  Hence there is a surjective holomorphic map $h:\C^n\to X$
  which has maximal rank somewhere, say at $0$
  (\cite{FF2017}).

  Since $Y$ is a Stein manifold, for a point $p\in Y$ we
  can find holomorphic functions $f_1,\ldots,f_n$ on $Y$
  such that $f_i(p)=0$ for all $i$ and
  \[
  \left(df_1\wedge\ldots\wedge df_n\right)_p\ne 0
  \]
  To find the desired map $F$ it suffices to take
  \[
  F=h\circ f,\quad f:q\mapsto (f_1(q),\ldots,f_n(q))
  \]
\end{proof}

\begin{proposition}\label{small-hom}
  Let $X$ be an $RR$-space, $E\subset X$ a finite subset and
  $p\in X\setminus E$.

  There exist $n\in\N$ and complete vector fields $w_j$ ($j\in\{1,\ldots, n\}$)
  on $X$ such that:
  \begin{enumerate}
  \item
    \[
    (w_j)_q=0\in T_qX\ \forall q\in E,\forall j\in\{1,\ldots, n\}
    \]
  \item
    The vector space $T_pX$ is generated by $(w_j)_p$
    ($j\in\{1,\ldots, n\}$)
  \end{enumerate}
\end{proposition}

\begin{proof}
  Follows from Proposition~\ref{t-span-vf} (ii) with $D=E\setminus\{p\}$.
\end{proof}

\begin{lemma}\label{step-dense}
  Let $X$ be an $RR$-space, $K$ a compact subset, $D\subset X$ a dense
  subset, $p\in X$ and $E \subset X \setminus\{p\}$ a finite subset.

  Assume that $X$ is equipped with a distance function.
  Let $\varepsilon>0$.
  
  Then there exists an automorphism $\phi$ of $X$ with
  \begin{enumerate}
  \item
    $\phi(q)=q\ \forall q\in E$,
  \item
    $\phi(p)\in D$.
  \item
    $d(x,\phi(x))<\varepsilon$  for all $x\in K$.
  \item
    $d(\phi^{-1}(x),x)<\varepsilon$  for all $x\in K$.
  \end{enumerate}
\end{lemma}

\begin{proof}
  Choose $n\in\N$ and complete vector fields $w_j$
  as in Proposition~\ref{small-hom}.
  Let $\phi^{(j)}_t$ denote the associated flows.
  Consider the map
  $\zeta:\C^ n\times X \to X$ defined as
  \[
  \zeta: (t_1,\ldots,t_n;x)\mapsto
  \left(\phi^{(1)}_{t_1}\circ\ldots\circ\phi^{(n)}_{t_n}\right)(x)
  \]
  By construction $\zeta(t;q)=q$ for all $t\in\C^ n$ and $q\in E$.
  Furthermore
  \[
  \C^ n\ni t\mapsto \zeta(t,p)\in X
  \]
  is submersive near the origin of $\C^ n$ since the tangent vectors
  $(w_j)_p$ span $T_pX$.

  Because $D$ is dense in $X$, we may choose a sequence $t^{(k)}\in\C^ n$
  with $\lim_{k\to\infty} t^{(k)}=0$ and $\zeta(t^{(k)};p)\in D$.
  Observe that the maps $x\mapsto  \zeta(t;x)$ converge locally uniformly
  to the identity map of $X$ for $t\mapsto 0$. Thus the condition
  \[
  d(x,\zeta(t,x))<\varepsilon\ \forall x\in K
  \]
  will be satisfied for $||t||$ small enough.

  Hence for suitably chosen $t$ the map $\phi:x\mapsto\zeta(t;x)$
  satisfies properties $(i)\sim(iii)$.

  The inverse of the automorphism
  $\left(\phi^{(1)}_{t_1}\circ\ldots\circ\phi^{(n)}_{t_n}\right)$
  is
  \[
  \left(\phi^{(n)}_{-t_n}\circ\ldots\circ\phi^{(1)}_{-t_1}\right)
  \]
  which also converges to the identiy map for $t\to 0$.
  Hence $(iv)$ is likewise satisfied for sufficiently small
  $t$.
\end{proof}

\begin{proposition}\label{countable}
  Let $S_1$ and $S_2$ be countable dense subsets of an $RR$-space  $X$.

  Then there exists an automorphism $\phi$ of $X$ with
  $\phi(S_1)=S_2$.
\end{proposition}

\begin{proof}
  Given Lemma \ref{step-dense}, we may use the same arguments
  as Rosay and Rudin in the proof of Theorem 2.2. of \cite{RR}.
\end{proof}
\section{Projection method.}

\subsection{Spreading families}

\begin{definition}\label{def-spread}
  Let $X$, $Y$ be locally compact topological spaces
  with exhaustion functions%
  \footnote{An exhaustion funtion is a continuous function
    $\rho:X\to\R^+_0$ such that $\rho^{-1}([0,c])$ is compact
    for every $c>0$.}
  denoted as $\rho_X$ resp.~$\rho_Y$.
  
  A {\em spreading family}
  is  a family
  of continuous mappings $\pi_t:X\to Y$ parametrized
  by elements $t$ in a probability space
  $(P,\mu)$
such that

\begin{equation}\label{star}
\forall r,\delta>0,\ \exists R>0,
\ \forall x\in X, \rho_X(x)>R:\quad
\mu\{t\in P: \rho_Y(\pi_t(x))<r\}<\delta
\end{equation}

\end{definition}

For simplyfing notation, in the sequel both $\rho_X$ and $\rho_Y$ are simply
denoted as $\rho$.

\begin{proposition}

  The above property \eqref{star} is equivalent to the
  following:

{\em
  For every $\delta>0$ and every compact subset $K$ in $Y$ there exists a
  compact subset $C$ in $X$ such that
  \[
  \forall x\in X\setminus C:
  \mu\left\{ \pi_t(x)\in K\right\}<\delta
  \]
}
\end{proposition}

\begin{proof}
$\implies$:

Choose $r>\sup\{\rho(p):p\in K\}$, choose $R$ as above and set
$C=\{p\in X:\rho(p)\le R\}$.
Then
\[
x\in X\setminus C\implies\rho(x)>R,\quad
\{t:\pi_t(x)\in K\}\subset\{t:\rho(\pi_t(x))<r\}
\]

$\Leftarrow$:

Given $r$, choose $K=\{y\in Y:\rho(y)\le r\}$.
Given $C\subset X$, let $R=\sup\{\rho(x):x\in C\}$.

Then
\[
\rho(x)>R\implies x\in X\setminus C,\quad
\{t:\rho(\pi_t(x))<r\}\subset \{t:\pi_t(x)\in K\}
\]
\end{proof}

\begin{corollary}
  The property of being as spreading family as in Definition~\ref{def-spread}
  is independent of the choice of the exhaustion functions $\rho_X$, $\rho_Y$
  used in this definition.
\end{corollary}

\begin{remark}
  Spreading families might be considered to be a probabilistic analog
  of proper maps. In the case of a trivial probability measure where
  $\mu(\{q\})=1$ for one point $q\in P$ we have a spreading family
  if and only if the map $\pi_q$ is a proper map.
\end{remark}

\begin{example}
  Let $X=\R^2$, $Y=\R$ and let $S^1$ act by rotations on $\R^2$.
  Let $\mu$ denote the normalized Haar measure on the compact
  group $S^ 1$.
  Let $p_1(x,y)=x$ and define for $P=(S^ 1,\mu)$
  \[
  \pi_t(v)=\left(p_1(t\cdot v)\right)
  \]

  Then $(\pi_t)$ is a {\em spreading} family
  (with the euclidean norm as exhaustion function
  $\rho$ on $X=\R^2$ resp.~$Y=\R$).

  In fact,
  \[
  |p_1\left(t\cdot v\right)|=|\cos(t)p_1(v)|
  \]
  implying
  \begin{align*}
  \mu\left\{ |\pi_t(v)|<r\right\}
  &= \mu\{ ||\cos(t)v||<r\}\\
&  = \mu\left\{ |\cos(t)| < \frac{r}{||v||} \right\}\\
&  = 4 \arcsin \left(\frac{r}{||v||} \right)\\
  \end{align*}
  Thus in order to fulfill \eqref{star} it suffices to
  choose $R$ such that
  \[
  \frac{r}{R} < \sin(\delta/4 )
  \]
\end{example}

\begin{proposition}\label{gen-proj}
  Let $X$, $Y$ be topological spaces with exhaustion functions $\rho$
  and a spreading family $\pi_t:X\to Y$ ($t\in (P,\mu)$).

  Then there exists a sequence of positive numbers $R_n$ with the following
  property:
  \begin{center}
  \fbox{
    \begin{minipage}{10 cm}
  {\em Let $p_n\in X$ be a sequence satisfying
    \[
    \forall n\in\N: \rho(p_{n+1})\ge\rho(p_n)>R_n
    \]
    and $D=\{p_n:n\in \N\}$.
    
    Then for almost all ($=$probability one) $t\in P$ the
    map $\pi_t$ maps $D$ properly ($=$ with finite fibers)
    onto a discrete
    subset of $Y$.
    }
    \end{minipage}
  }
  \end{center}
\end{proposition}

\begin{proof}
  For every $n\in\N$ we set
  $\delta_n=2^{n+1}$. Since $\pi_t$ is a spreading family,
  there is a number $R_n>0$ such that
  \[
  \forall x\in X, \rho(x)>R_n:\quad
  \mu\{t\in P: \rho(\pi_t(x))<n)\}<\delta_n
  \]
  (Definition~\ref{def-spread} with $r=n$.)

If $\pi_t(D)$ is not discrete, 
there exists a radius $r>0$ such that $\rho(\pi_t(p_n))<r$
for infinitely many $n\in\N$. But for every $\epsilon>0$ there
exists a number $N$ such that
\[
\sum_{n\ge N}\delta_n
=\sum_{n\ge N}2^{-(n+1)}=2^{-N}<\epsilon,
\]
which implies that
\[
\mu\{t\in P: \exists n>N: \rho(\pi_t(p_n))>n\}<\epsilon.
\]
This in turn implies that
\[
\mu\{t\in P:\# (B_r\cap\pi_t(D))>N\}<\epsilon.
\]
Since this holds for all $\epsilon>0$, we know that
$\pi_t^{-1}(B_r)\cap D$
is infinite with probability $0$. Using the $\sigma$-additivity
of the measure $\mu$, it follows that
\[
\cup_{r\in\N}\left\{ t\in P: \# (B_r\cap\pi_t(D))=\infty\right\}
\]
is a set of measure zero.
Hence $\pi_t$ maps $D$ onto a discrete set in $Y$ with finite
fibers (i.e.,$\pi_t|_D$ is proper) with probability one.
\end{proof}

\subsection{Threshold Sequences}

As in \cite{JW-TAME1} and \cite{JW-LAG},
we use the following notion:

\begin{definition}\label{thres-def}
  Let $X$ be a Stein manifold with an exhaustion function $\rho$.
  An increasing sequence $R_n$ of positive numbers is called
  {\em threshold sequence} iff every discrete subset $D\subset X$
  satisfying
  \[
  \forall n\in\N: \#\{ x\in D:\rho(x)\le R_n\} \le n
  \]
  is a (weakly) tame set.
\end{definition}

\begin{theorem}\label{spread-thres}
Let $X$, $Y$ be complex manifolds, $Y$ Stein, $\pi_t:X\to Y$ be a family
of holomorphic mappings, parametrized by a probability space
$t\in (P,\mu)$ which is a {\em spreading family} in the sense of
Definition~\ref{def-spread}.

Assume that for every $t$ there is a complex Lie group
$G_t$ with a holomorphic group action
$G_t\times X\to X$ such that every $G_t$-orbit in $X$
is a non-compact
closed subset in a $\pi_t$-fiber.

Then there exists a {\em threshold sequence}.
\end{theorem}

\begin{proof}
  Fix an exhaustion function $\rho$ on $X$.
  Choose $R_n$ as in Proposition~\ref{gen-proj}.
  Let $D$ be a discrete subset satisfying
  \[
  \forall n\in\N: \#\{ x\in D:\rho(x)\le R_n\} \le n.
  \]

  Fix a projection $\pi_t:X\to Y$ such that
  $\pi_t|_D$ is injective with discrete image.

  Let $\zeta:D\to\R^ +$ be an arbitrary map.
    
  Since every $G_t$-orbit is unbounded,
  for every $\gamma\in D'=\pi_t(D)$ there is an
  element
  $g_\gamma\in G_t$ for 
  such that
  \[
  \rho\left(g_\gamma(\tilde\gamma)\right)>\zeta(\tilde\gamma)
  \quad\forall \tilde\gamma\in D\text{ with }\pi_t(\tilde\gamma)=\gamma
  \]
  (Here we use the fact that
  $\{\tilde\gamma\in D:\pi_t(\tilde\gamma)=\gamma$ is finite.)

  We obtain a map $\xi:D'\to G_t$ defined as $\xi(\gamma)=g_\gamma$.
  Because $Y$ is Stein and $G_t$ is a complex Lie group, this map
  extends to a holomorphic map $\eta:Y\to G_t$. The holomorphic
  automorphism $\phi$ of $X$ defined as
  \[
  \phi(x)=\eta\left(\pi_t(x)\right)(x)
  \]
  satisfies
  \[
  \rho(\phi(x))\ge \zeta(x)\quad\forall x\in D
  \]
  Thus $D$ is tame and $(R_n)_n$ is a threshold sequence.
\end{proof}

\subsection{Flexible Varieties}

For basic facts on {\em flexible} affine varieties, see \cite{AFKKZ}.
Flexible varieties are rich in ``LNDs''. These are {\em locally nilpotent
derivations} of the function ring. They correspond to complete
holomorphic vector fields whose flow defines an algebraic action of the additive
group $(\C,+)$ on the given variety. By abuse of language we call these
vector fields also ``LNDs''.

\begin{proposition}\label{flex-tame}
  Let $X$ be a smooth complex flexible affine variety, $C\subset X$ an
  algebraic curve isomorphic to the affine line $\C$ and $D\subset C$ an infinite
  discrete subset.

  Then $D$ is a weakly tame discrete subset of $X$.
\end{proposition}

\begin{proof}
  By definition, the tangent space of a smooth flexible variety is everywhere
  spanned by LNDs.
  Fix $p\in C$ and let $v_1,\ldots,v_n$ be LNDs spanning $T_pX$.
  Then $v_1,\ldots,v_n$ span $TX$ at  every point in a Zariski open
  subset $C^*$ of $C$ (where ``Zariski open'' refers to the algebraic Zariski
  topology). Thus $C\setminus C^*$ is finite and it is clear that we may
  choose finitely many LNDs $v_1,\ldots,v_N$ which span $TX$ at every point
  of $C$.
  Using the flows $\phi_{j,t}$ associated to $v_j$ we obtain
  a morphism $\Psi:\C^N\times X\to X$ defined as
  \[
  \Psi:(t_1,\ldots,t_n;x)\mapsto
  \left(\phi_{1,t_1}\circ\ldots\circ\phi_{N,t_n}\right)(x)
  \]
  By construction, this morphism is submersive along $\{0\}\times C$.
  
  Next we fix one LND $w$ on $X$.
  The invariants for any non-trivial algebraic action of
  the additive group $(\C,+)$ on $X$
  generically separate orbits due to the theorem of Rosenlicht,
  see \cite{PopovVinberg}, Theorem~2.3.
  In particular (since $\dim(X)\ge 2$ because $X$ is flexible)
  there exists a non-constant polynomial function $P$ which is invariant
  under the flow of $w$, i.e., $w(P)\equiv 0$.

  Let $Z$ denote the vanishing locus of $w$. For $q\in D$ we define
  \[
  W_q=\{t\in\C^N:\Psi(t,q)\not\in Z\}
  \]
  and
  \[
  \Omega=\{t\in\C^N: \text{$P\circ\Psi(t,-)$ is not constant
  on $C$}\}\cap\left( \bigcap_{q\in D}W_q\right)
  \]
  Note that $\Omega$ is the complement of a union of countably many
  strict subvarieties of $X$ and therefore not empty.

  By replacing $w$ with $\Psi_t^*w$ for s suitably choosen $t$, we
  may thus assume that
  \begin{itemize}
  \item
    $P$ is not constant on $C$.
  \item
    $w$ does not vanish at any $q\in D$.
  \end{itemize}

  Fix an exhaustion function $\rho$ on $X$ and a function
  $\zeta:D\to\R^+$.
  In order to verify that $D$ is weakly tame, we have to show that
  there exists a holomorphic automorphism $\alpha$ of $X$ such that
  \[
  \rho(\alpha(q))>\zeta(q)\ \forall q\in D
  \]

  The restriction of $P$ to $C$ defined a non-constant polynomial function
  from $C\simeq\C$ to $\C$ and is therefore proper.
  Thus $D'\dfe P(D)$ is a discrete subset of $\C$ and
  the fibers of $P|_D:D\to D'$ are finite.

  Consider the flow $\phi_t$ associated
  the vector field $w$. Recall that $w$ does not vanish at any point of $D$.
  Thus every orbit $\omega$ of the flow through a point of $D$ is isomorphic
  to an affine line and in particular it is a closed subset of $X$.
  Hence
  \[
  \{x\in\omega:\rho(x)\le c\}
  \]
  is compact for every $c\in\R$ and every orbit $\omega$.
  
  It follows that
  \[
  K_q=\{ t\in\C: \rho(\phi_t(q))\le\zeta(q)\}
  \]
  is compact for every $q\in D$.

  Since $P:D\to D'$ has finite fibers,
  \[
  \cup_{q\in P^{-1}(z)}K_q
  \]
  is compact for every $z\in D'$.
  Hence for every $z\in D'$ we may choose a complex number
  $s(z)$ such that
  \[
  \forall q\in P^{-1}(z):s(z)\not\in K_q
  \]
  This in turn implies
  \[
  \forall q\in D:
  \rho\left(\phi_{s(P(q))}(q)\right)>\zeta(q)
  \]

  Since $D'$ is discrete in $\C$ we may find a holomorphic function
  $h:\C\to\C$ with $h(z)=s(z)\ \forall z\in D'$.
  Now $h\circ P$ is an $w$-invariant function, which implies that
  $\tilde w=(h\circ P)\cdot w$ is a {\em complete}
  holomorphic vector field.

  By construction, the time $1$-map $\alpha$ of the flow associated
  to $\tilde w$ satisfies
  \[
  \rho(\alpha(q))>\zeta(q)\ \forall q\in D
  \]
  
\end{proof}

\begin{corollary}
  Every flexible smooth complex affine variety contains a weakly tame discrete
  subset.
\end{corollary}

\begin{proof}
  Let $X$ be a flexible complex affine variety. Any non-trivial orbit of
  an algebraic $(\C,+)$-action is closed and isomorphic to $\C$.
  Thus $X$ contains an algebraic curve $C$ which is isomorphic to the
  affine line and we may use Proposiiton~\ref{flex-tame}.
\end{proof}

\begin{theorem}\label{thm-flex-tame}
  Let $X$ be a flexible smooth affine variety which is an $RR$-space
  as complex manifold.

  Let $C\subset X$ be an algebraic curve isomorphic to $\C$ and
  let $D_0\subset C$ be an infinite discrete subset.

  Then for any infinite discrete subset $D\subset X$ the following
  conditions are equivalent:
  \begin{enumerate}
  \item
    $D$ is tame.
  \item
    There is a holomorphic automorphism of $X$ with $\phi(D)\subset C$.
  \item
    There is a holomorphic automorphism of $X$ with $\phi(D)=D_0$.
  \end{enumerate}
\end{theorem}

\begin{proof}
  $(iii)\implies(ii)$ is trivial.
  $(ii)\implies (i)$ follows from Proposition~\ref{flex-tame}.

  $(i)\implies(iii)$.
  First we observe that $D_0$ is tame due to Proposition~\ref{flex-tame}.
  Since $X$ is an $RR$-space,  $(i)\implies(iii)$ now follows
  from Axiom $(RR1)$.  
\end{proof}

\begin{corollary}
  Let $P\in\C[X]$ with $P(0)=0$ and without multiple roots.
  Let $X=\{(x,y,z)\in\C: xy=P(z)\}$
  and $D_0=\{(n,0,0):n\in\N\}$

  Then an infinite discrete subset $D\subset X$ is tame if and only if
  there is a holomorphic automorphism $\phi$ of $X$ with $\phi(D)=D_0$.
\end{corollary}

\begin{proof}
  Note that
  \[
  D_0\subset C=\{(x,0,0):x\in\C\}\subset X
  \]
  and that $X$ is both flexible and an $RR$-space.
\end{proof}

\section{Danielewski surfaces}

Let $P$ be a polynomial without multiple roots.
Then
\[
X=\{(x,y,z)\in\C:xy=P(z)\}
\]
is called a {\em Danielewski surface}.

  Danielewski surfaces are
  known to be flexible affine varieties and to satisfy the density
  property (\cite{KalimanKutz2008.1}). One may easily verify that
  $b_2(X)=\deg(P)-1$ for a Danielewski surface
  (compare \cite{MR1669174}).
  
We summarize some key properties.

$X$ is smooth (because $P$ has no multiple roots),
simply-connected (if $\deg(P)>0)$) with
$\chi(X)=\deg(P)$ and $b_2(X)=\deg(P)-1$.

If $\deg(P)=1$, we have $X\simeq\C^2$. For $\deg(P)=2$ the surface $X$ is a
smooth affine quadric on which $SL_2(\C)$ acts transitively.
However, for $\deg(P)>2$ the automorphism group acts transitively
on $X$ although no finite-dimensional Lie group does:

  The classification
  of homogeneous surfaces
  (see \cite{MR638369}) implies that no simply-connected
  complex surface with $b_2(X)>1$ is biholomorphic to a quotient $G/H$ of
  complex Lie groups.

Now $\theta_1\dfe P'(z)\partial_y-x\partial_z$ is a LND with
$x$ as invariant, and $\theta_2\dfe P'(z)\partial_x-y\partial_z$
is a LND with $y$ as invariant.

Observe that $(x,y,z)\mapsto(y,x,z)$ ist an automorphism of $X$
which interchanges $\theta_1$ and $\theta_2$.

The LND $P'(z)\partial_x-y\partial_z$ integrates to the $1$-parameter group
\begin{align*}
(x,y,z)\mapsto &\left( x+\sum_{k=1}^\infty (-y)^{k-1}t^k\frac{P^{(k)}(z)}{k!},
y, z-yt\right)\\
&=
\begin{cases}\left( \frac 1yP(z-yt),y,z-yt\right) & \text{ if $y\ne 0$}\\
  (x+tP'(z),0,z) & \text{ if $y=0$}\\
\end{cases}\\
\end{align*}
(The equality of the two forms follows from $xy=P(z)$.)

If $x\ne 0$, the level set $x=c$ agrees with one orbit
associated to the flow of $\theta_1$:
\[
\left\{(x,y,z):z\in\C,x=c,y=\frac{P(z)}{c}\right\}\simeq\C.
\]

The level set $x=0$  is a union of $\deg(P)$-many orbits:
\[
\{(0,y,z):y\in\C, P(z)=0\}
\]

For every polynomial $Q\in\C[X]$ we have the ``replicas''
\[
Q(x)\left( P'(z)\partial_y-x\partial_z \right),
\quad
Q(y)\left( P'(z)\partial_x-y\partial_z \right),
\]
serving as analogs for shears.

Non-unipotent automorphisms can be constructed via
\[
(x,y,z)\mapsto (xe^{\phi(z)},ye^{-\phi(z)},z)
\]
for any holomorphic function $\phi$.

Summarizing:

\begin{proposition}\label{dan-action}
  Let $X$ be a Danielewski surface. Then there
  exist regular functions $p_i:X\to \C$ ($i=1,2$) and actions
  $\mu_i: G\times X\to X$ ($i=1,2$)
  of the additive group $G=(\C,+)$ such that:
  \begin{enumerate}
  \item
    $p_i$ is a $\mu_i$-invariant, i.e.,
    $p_i(\mu_i(t,x))=p_i(x)$
  \item
    For all $t\in\C\setminus\{0\}$, $i=1,2$ the fiber
    $p_i^{-1}(t)$ is one $(G,\mu_i)$-orbit.
  \item
    Each $p_i$ maps each $G_{3-i}$-orbit onto $\C$.
    In particular, every orbit is one-dimensional; there are no fixed points.
  \item $x\mapsto (p_1(x),p_2(x))$
    defines a finite morphism from $X$ to $\C^2$.
  \end{enumerate}
\end{proposition}

\begin{proposition}\label{dan-spread}
Let $\C$ be endowed with the Gaussian measure $\gamma$, i.e., 
\[
\gamma(U)=\int_U
\frac1{2\pi }e^{-\frac{r^2}{2}}dxdy\quad
\text{ for $U\subset\C$, $r=|z|=\sqrt{x^2+y^2}$}.
\]
Let $P$ be a polynomial of degree at least two without multiple
roots and let
\[
X=\{(x,y,z)\in\C^3:xy=P(z)\}
\]
be the associated Danielewski surface.

Define 
\[
\eta_t(x,y,z)\mapsto
\begin{cases}
\frac 1yP(z+yt) & \text{if $y\ne 0$}\\
x+tP'(z) & \text{if $y= 0$}\\
\end{cases}
\]

Then for every $r>0$, $\epsilon>0$, there exists a number $R>0$
such that for 
each $(x,y,z)\in X$ with $\max\{|x|,|y|\}\ge R$ we have
\[
\gamma\left( \{t\in\C: 
|\eta_t(x,y,z)|<r\} \right) < \epsilon
\]
\end{proposition}

\begin{lemma}\label{eta-non-constant}
  With the assumptions and notations of Proposition~\ref{dan-spread},
  for every $(x,y,z)\in X$ the map $t\mapsto\eta_t(x,y,z)$
  is a non-constant
  function.
\end{lemma}

\begin{proof}
  If $y\ne 0$, we have $\eta_t(x,y,z)=\frac 1yP(z+yt)$ and the assertion
  holds, because $P$ is not constant and $y\ne 0$.

  If $y=0$, we have $\eta_t(x,y,z)=x+tP'(z)$. In this case
  $t\mapsto \eta_t(x,y,z)$ is constant iff $P'(z)=0$.
  Now $P$ has no zeroes with multiplicity $\ge 2$, hence
  $P'(z)=0\ \implies\ P(z)\ne 0$. But
  \[
  (x,y,z)\in X\ \iff xy=P(z)
  \]
  Hence $P(z)\ne 0\ \implies\ y\ne 0$.
  Therefore the case $y=0= P'(z)$ is not possible.
\end{proof}

\begin{remark}
  For a Danielewski surface
  \[
  X=\{(x,y,z):xy=P(z)\}
  \]
  the map $(x,y,z)\mapsto\max\{|x|,|y|\}$ is an exhaustion function on $X$.
  Indeed, if $|x|,|y|\le R$, then $|P(z)|=|xy|\le R^2$ and for every
  polynomial the set $\{z:|P(z|\le r\}$ is compact for all $r$.

  Hence the above Proposition~\ref{dan-spread}
  states that $\eta_t$ is a {\em spreading family}
  for $t\in (\C,\gamma)$.
\end{remark}

\begin{proof}[Proof of Proposiiton~\ref{dan-spread}]
  We need an estimate for the polynomial  $P$.
  Let $d$ be the degree of $P$.
  There are constants
$\rho,M,\alpha,\beta>0$ such that
\[
\alpha|z|^d\le |P(z)|\le \beta|z|^d\quad
\text{ if $|z|\ge\rho$}
\]
and
\[
|P(z)|\le M \quad \text{ if $|z|\le\rho$}.
\]

We claim:

\begin{claim}\label{claim-1}
  
  \begin{center}
  \fbox{
    \begin{minipage}{10 cm}
{\em For every $(x,y,z)\in X$ and $r>0$ we have
  \begin{equation*}
  \gamma\left\{ t:\left|\eta_t(x,y,z)\right|<r\right\}
  \le  \pi\frac{\rho^2}{|y|^2}
  +|y|^{\frac{2-2d}{d}}\pi
  \left(\frac r\alpha\right)^{\frac 2d}
  \end{equation*}
  }
    \end{minipage}
  }
  \end{center}
\end{claim}

We observe that
the Lebesgue measure $\lambda$  is an upper bound
for the Gaussian measure $\gamma$, because 
$\frac{1}{2\pi}e^{-(x^2+y^2)/2}< 1$
for every $(x,y)\in\R^2$.

Hence
\[
\gamma\{|z+yt|\le\rho\}
\le
\lambda\{|z+yt|\le\rho\}
= \pi\frac{\rho^2}{|y|^2}
\]
and
\begin{align*}
& \gamma\left\{t:|z+yt|>\rho,\ \left|\frac{1}{y}P(z+yt)\right|<r\right\}\\
\le & \gamma\left\{t:|z+yt|>\rho,\ \alpha\left|\frac{1}{y}(z+yt)^d\right|<r\right\}\\
\le & \gamma\left\{t :\alpha\left|\frac{1}{y}(z+yt)^d\right|<r\right\}\\
& \text{because $|z+ty|\ge \rho\implies |P(z+yt)|\ge \alpha |z+yt|^d$}\\
\le & \lambda\left\{t :\alpha\left|\frac{1}{y}(z+yt)^d\right|<r\right\}\\\le & \lambda\left\{t :\left|(z+yt)^d\right|< \frac{|y|r}{\alpha} \right\}
  = \lambda\left\{t :|(z+yt)|< \left|\frac{yr}{\alpha}\right|^{\frac 1d} \right\} \\
  = & \lambda\left\{t :\left|\frac zy +t\right|<|y|^{\frac 1d-1}
  \left|\frac{r}{\alpha}\right|^{\frac 1d}\right\} \\ 
= & \pi \left ( \frac{r}{\alpha} \right)^{\frac 2d}|y|^{\frac 2d-2}
      \\
\end{align*}
Hence
\begin{align*}
&\gamma\left\{ t:\left|\eta_t(x,y,z)\right|<r\right\}\\
&=
\gamma\left\{t:\left|\frac{1}{y}P(z+yt)\right|<r\right\}\\
&\le \gamma\{|z+yt|\le\rho\}+
\gamma\left\{t:|z+yt|>\rho,\ \left|\frac{1}{y}P(z+yt)\right|<r\right\}\\
&\le  \pi\frac{\rho^2}{|y|^2}
+ \pi \left ( \frac{r}{\alpha} \right)^{\frac 2d}|y|^{\frac 2d-2} \\
\\
\end{align*}

This proves the first claim.

Next we fix a positive number $r>0$.
For $R>0$  we define $c_R$ as
\[
c_R=\inf\{ |t| :\exists (x,y,z)\in X, |\eta_t(x,y,z)|\le r,
|x|=R, |y|\le R^{1/d^2}
\}
\]
We claim:
\begin{claim}\label{claim-2}
\[
\lim_{R\to\infty}c_R=+\infty
\]
\end{claim}
To verify this, assume the contrary.

Then we have sequences $t_n$, $R_n$, $x_n$, $y_n$, $z_n$
such that
\begin{enumerate}
\item
  $\lim R_n=+\infty$,
\item
$C=\sup_n|t_n|<\infty$.
\item
$(x_n,y_n,z_n)\in X$,
\item
$|x_n|=R_n$,
\item
$|y_n|\le R_n^{1/d^2}$,
\item
$|\eta_{t_n}(x_n,y_n,z_n)|\le r$.
\end{enumerate}

Define
\[
Q(z,y,t)=\frac{1}{y}\left( P(z)-P(z+yt)\right)
\]
Recall $d=\deg(P)\ge 2$. We note that
$Q$ is a polynomial in degree $d\ge 2$ in $t$, and in degree $d-1\ge 1$ in
$y$ and $z$.
Hence $\exists K>0$:
\[
|Q(z,y,t)| \le K (|z|+1)^{d-1}(|y|+1)^{d-1}(|t|+1)^d
,\ \ \forall z,y,t\in\C
\]
Now $|x_n|=R_n$ in combination with $(x_n,y_n,z_n)\in X$ implies
\[ 
R_n|y_n|=|P(z_n)|
\]
which in turn implies
\[
\left| \frac{1}{y_n}P(z_n)\right|=R_n.
\]
On the other hand
\[
\left|\frac 1 {y_n}P(z_n+y_nt_n)\right|\le r
\]
It follows that
\[
|Q(z_n,y_n,t_n)|\ge R_n-r
\]
But
\begin{align*}
|Q(z_n,y_n,t_n)| &\le K (|z_n|+1)^{d-1}(|y_n|+1)^{d-1}(|t_n|+1)^d \\
&\le \tilde K (|z_n|+1)^{d-1} \left(\left(R_n\right)^{\frac 1 {d^2}}+1\right)^{d-1}\\
\end{align*}
for $\tilde K=K(C+1)^d$
(using $|y_n|\le R^{1/d^2}$ and $|t_n|\le C$).

Since $P$ is a polynomial of degree $d$ we have
\[
\frac{1}{a_d}=\lim_{z\to\infty}\frac{(|z|+1)^d}{|P(z)|},
\]
where $a_d$ is the leading coefficient of the polynomial
$P(z)=\sum_k a_kz^k$.

It follows that for every $c'>\frac{1}{a_d}$ there is a constant
$N>0$ such that
\begin{align*}
&(|z_n|+1)^{d} \le c' |P(z_n)|&\forall n\ge N\\
\implies &(|z_n|+1)^{d-1} \le c' |P(z_n)|^{\frac{d-1}{d}}&\forall n\ge N\\
\end{align*}
Therefore
\[
|Q(z_n,y_n,t_n)|
\le
\tilde K c' |P(z_n)|^{\frac{d-1}{d}} \left(\left(R_n\right)^{\frac 1 {d^2}}+1\right)^{d-1}
 \]
Combined with $|P(z_n)|=|x_ny_n|\le R_n\cdot R_n^{\frac{1}{d^2}}
  =R_n^{1+\frac{1}{d^2}}$
we obtain
\[
|Q(z_n,y_n,t_n)|
\le
\tilde K c' \left|R_n^{1+\frac{1}{d^2}}\right|^{\frac{d-1}{d}}
    \left(R_n^{\frac 1 {d^2}}+1\right)^{d-1}
 \]
 Since $\lim R_n=+\infty$, there is no loss in generality in
 assuming
 \[R_n\ge 1\forall n
 .
 \]
 Then $R_n+1\le 2R_n$ and we
 may deduce
 \[
|Q(z_n,y_n,t_n)|
\le
 2^{d-1}\tilde K c' \left|R_n^{1+\frac{1}{d^2}}\right|^{\frac{d-1}{d}}
 \left(R_n^{\frac 1 {d^2}}\right)^{d-1}
 \]
 An easy calculation shows
 \[
 \left(1+\frac{1}{d^2}\right)\frac{d-1}{d}+\frac{d-1}{d^2}=1-\frac{1}{d^3}
 \]
 Thus
 \[
|Q(z_n,y_n,t_n)|
\le
 2^{d-1}\tilde K c' \left(R_n\right)^{1-\frac 1{d^3}}\quad\forall n\ge N
 \]
 On the other hand $|Q(z_n,y_n,t_n)|\ge R_n-r$.

 But
 \[
 R_n-r
\le
 2^{d-1}\tilde K c' \left(R_n\right)^{1-\frac 1{d^3}}\quad\forall n\ge N
 \]
 leads to a contradiction for $n\to\infty$,
 because $\lim_{n\to\infty}R_n=\infty$ and $1-\frac{1}{d^3}<1$.
 
Therefore $\lim_{R\to\infty}c_R=+\infty$, i.e., we proved claim~\ref{claim-2}.

\begin{claim}\label{claim-3}
  For every $\varepsilon>0$ there exists $R_0>0$ with the following property:
    \begin{center}
  \fbox{
    \begin{minipage}{10 cm}
      {\em For every $(x,y,z)\in X$ with
        $|x|\ge R_0>|y|^{d^2}$ we have
  \begin{equation*}
  \gamma\left\{ t:\left|\eta_t(x,y,z)\right|<r\right\}
  \le  \varepsilon
  \end{equation*}
  }
    \end{minipage}
  }
  \end{center}
\end{claim}

Let us prove this claim. Suppose $\varepsilon>0$ is given.

Recall that
\[
c_R=\inf\{ |t| :\exists (x,y,z)\in X, |\eta_t(x,y,z)|\le r,
|x|=R, |y|\le R^{1/d^2}
\}
\]

Thus
\[
\{ t:\exists (x,y,z)\in X, |\eta_t(x,y,z)|\le r,
|x|=R, |y|\le R^{1/d^2}
\}\subset
\{
t\in\C:|t|\ge c_R\}
\]

Observe that
\[
\gamma\left( \{t\in\C:|t|\ge s \}\right)
=\int_s ^\infty re^{-r^2/2}dr= 1 - e^{-s ^2/2}
\]
which implies
\[
\lim_{s \to\infty}\gamma\left( \{t\in\C: |t|\ge s \}\right)
=0
\]

We choose $R_0>0$ such that
\[
\forall R\ge R_0:\ \gamma
\left(\left \{t\in\C: |t|\ge c_R \right\}\right) <\varepsilon
\]
which is possible, because $\lim_{R\to\infty}c_R=0$.

As a consequence, we obtain:

{\em For every $R\ge R_0$ and every $(x,y,z)\in  X$ with
  $|x|=R$ and $|y|\le R^{1/d^2}$ we have
\[
 \gamma
\left(\left\{ t: |\eta_t(x,y,z)|\le r,
\right\}\right)<\varepsilon
\]
}

Therefore:

\[
\forall (x,y,z)\in X:
|x|\ge R_0\ge |y|^{d^2}
\implies
 \gamma
\left(\left\{ t: |\eta_t(x,y,z)|\le r,
\right\}\right)<\varepsilon
\]

On the other hand (Claim~\ref{claim-1}):
\[
  \gamma\left\{ t:\left|\eta_t(x,y,z)\right|<r\right\}
  \le  \pi\frac{\rho^2}{|y|^2}
  +|y|^{\frac{2-2d}{d}}\pi
  \left(\frac r\alpha\right)^{\frac 2d}
  \]

  Thus, in order to prove the assertion of the proposition
  it suffices to show  that
  \[
    \pi\frac{\rho^2}{|y|^2}
  +|y|^{\frac{2-2d}{d}}\pi
  \left(\frac r\alpha\right)^{\frac 2d}
  \le\varepsilon
  \]
  for all
  \[
  (x,y,z)\in X\setminus \{(x,y,z):|x|\ge R_0\ge |y|^{d^2}\}
  \]
  \begin{claim}\label{claim-4}
    Assume $R_0\ge 1$, $\max\{|x|,|y|\}\ge R_0$
    and $(x,y,z)\in X$.

    Then (at least) one of the following conditions hold:
    \begin{enumerate}
    \item
      $|x|\ge R_0\ge |y|^{d^2}$
    \item
      $|y|^{d^2}\ge R_0$.
    \end{enumerate}
  \end{claim}

  Indeed, it the first condition fails, $x<R_0$ or $R_0<|y|^{d^2}$
  holds. If $R_0<|y|^{d^2}$, we have condition $(ii)$.
  Hence we assume $x<R_0$.
  Now $x<R_0$ in combination with
  $\max\{|x|,|y|\}\ge R_0$ implies $|y|\ge R_0$.
  Next we observe that $|y|^{d^2}\ge|y|$, because of
  $|y|\ge R_0\ge 1$.
  Hence $|y|^{d^2}\ge R_0$ and again we have condition $(ii)$.
  This prove our claim~\ref{claim-4}.

Observe that
\[
\lim_{|y|\to\infty}\left(    \pi\frac{\rho^2}{|y|^2}
  +|y|^{\frac{2-2d}{d}}\pi
  \left(\frac r\alpha\right)^{\frac 2d}\right)=0,
  \]
because $d>1\implies \frac{2-2d}{d}<0$.

Hence, if we choose $R_0$ sufficiently large,
we obtain
\[
 \gamma\left\{ t:\left|\eta_t(x,y,z)\right|<r\right\}<\varepsilon\ge 1
\]
for every $(x,y,z)\in X$ with $\max\{|x|,|y|\}\ge R_0$.
\end{proof}

\begin{theorem}\label{dan-thres}
  Every Danielewski surface
  \[
  X=\{(x,y,z)\in\C^3:xy=P(z)\}
  \]
  admits a {\em threshold sequence}.
\end{theorem}

\begin{proof}
Define 
\[
\eta_t(x,y,z)\mapsto
\begin{cases}
\frac 1yP(z+yt) & \text{if $y\ne 0$}\\
x+tP'(z) & \text{if $y= 0$}\\
\end{cases}
\]
Due to Proposition~\ref{dan-spread} the family of maps
$\eta_t:X\to\C$ parametrized by $\C$, equipped with the Gaussian
measure is a {\em spreading family} in the sense of
Definition~\ref{def-spread}.

Thus theorem~\ref{spread-thres} implies that $X$ admits a threshold
sequence.
\end{proof}

\begin{proposition}\label{tame-proper-proj}
  Let $D$ be a weakly tame discrete subset in a Danielew\-ski
  surface
  \[
  X=\{(x,y,z)\in\C^3:xy=P(z)\}
  \]
  Let $\pi:(x,y,z)\mapsto x$ be the natural projection.
  
  Then there exists a holomorphic automorphism $\phi$ of $X$
  such that $\pi\circ\phi$ restricts to a proper map from $D$
  to $\C$ which avoids $0$ (i.e.~$0\not\in\pi(\phi(D))$).
\end{proposition}

\begin{proof}
  As verified in the proof of the above theorem~\ref{dan-thres}
  $\eta_t$ is a spreading family. By construction
  $\eta_t=\pi\circ\phi_t$ for some holomorphic automorphism
  $\phi_t$ of $X$.

  Fix an exhaustion function $\rho$ on $X$.
  Due to Proposition~\ref{gen-proj} there are constants $R_n$ such that
  there exists a $t$ with $\eta_t|_D$ proper if
  \begin{equation}\label{rnx}
    \#\{p\in D:\rho(p)\le R_n\}\le n\ \forall n
  \end{equation}
  $(1)$. Due to the definition of the notion {\em ``weakly tame''} there
  is a holomorphic automorphism $\psi$ of $X$ such that
  $\psi(D)$ satisfies the condition \eqref{rnx}.
  It follows that for  a randomly chosen $t$ the map $\eta_t\circ\psi$
  restricts to a proper map on $D$ with probability $1$.
  (Proposition~\ref{gen-proj})
  
  $(2)$. Fix an element $(x,y,z)\in \psi(D)$.
  Due to Lemma~\ref{eta-non-constant}
  the set
  \[
  \{t\in\C: \eta_t(x,y,z)=0\}
  \]
  is finite and therefore has measure zero.
  Since $D$ is countable, it follows that for a randomly chosen
  parameter $t$ we have
  \[
  0\not\in\eta_t(\psi(D))
  \]
  with probability $1$.
  
  $(3)$ Combined, these two observations yield the assertion
  of the proposition, since $\phi_t\circ\psi\in\Aut(X)$.
\end{proof}

\begin{proposition}
  Let $D$ be a weakly tame discrete subset in a Danielew\-ski
  surface
  \[
  X=\{(x,y,z)\in\C^3:xy=P(z)\}
  \]

  Let
  \[
  p_1:(x,y,z)\mapsto y,\quad p_2:(x,y,z)\mapsto x
  \]
  be the natural projections.

  Then there exists a holomorphic automorphism $\phi$ of $X$ such that
  $p_1\circ\phi$ restricts to a proper and injective map
  from $D$ to $\C$.
\end{proposition}

\begin{proof}
  There is an automorphism $\psi$ such that $p_2\circ\psi:D\to\C$ is
  proper and avoids $0$ (Proposition~\ref{tame-proper-proj}).
  Thus we may without loss of generality assume that $p_2|_D:D\to\C$
  is proper and $p_2(D)\subset\C^*$.
  
  Enumerate $p_2(D)=\{q_k:k\in\N\}$.
  Since $p_2$ restricts to a proper map on $D$, the set
  \[
  \{(x,y,z)\in D:x=q_k,\}
  \]
  is finite for every $k\in\N$.
  
  Let $\mu_t$ denote the action along the $p_2$-fibers as given by
  \[
  \mu_t:(x,y,z)\mapsto  \left(x,\frac{1}{x}P(z-xt),z-xt\right)
  \]
  From Proposition~\ref{dan-action} and $q_k\ne 0\forall k$
  we may deduce that $p_2(q_k)$ is a single orbit
  of $\mu$ for every $k\in\N$.
  
  Recursively we choose numbers $c_k\in\C$, $r_k\in\R^+$ such that
  \begin{enumerate}
  \item
    $p_1\circ \mu_{c_k}$ is injective on
    $p_2^{-1}(q_k)\cap D$.
  \item
    \[
    \left|p_1(v)\right|>r_{k-1}\quad\forall v\in D\cap p_2^{-1}(q_k)
    \]
  \item
    \[
    r_k=\max\left\{|p_1(v)|:v\in p_2^{-1}\left(q_1,\ldots q_k\right)\right\}
    \]
  \end{enumerate}
  
  We choose a holomorphic function $\zeta:\C\to\C$ such that
  $\zeta(q_k)=c_k$ and we define $\alpha$ to be the time-one-map
  of the holomorphic flow associated to $(\zeta\circ p_2)\delta$
  where $\delta$ denotes the LND inducing the flow $\mu_t$, i.e.,
  $\delta=P'(z)\partial_y-x\partial_z$.
  Note that $(\zeta\circ p_2)\delta$ is a complete holomorphic vector
  field, because $\delta(x)=0$ and consequently
  $\delta(\zeta\circ p_2)=0$.
  Now $\alpha(q)=\mu_{c_k}(q)$ for every $q$ with $p_2(q)=q_k$
  for every $k\in\N$. 
  
  Then $p_1\circ \alpha\circ\psi$ maps $D$ injectively
  and properly to $\C$.
\end{proof}

\begin{theorem}
  Let
  \[
  S=\{xy=P(z)\}
  \]
  be a Danielewski surface with an infinite discrete
  subset $D$. Let $p:S\to\C$ be the projection morphism
  $p:(x,y,z)\mapsto x$.

  Then the following properties are equivalent:
  \begin{enumerate}
  \item
    $D$ is weakly tame.
  \item
    $D$ is strongly tame.
  \item
    There is a holomorphic automorphism $\phi$ of $S$ such that
    $p\circ\phi:D\to \C$ is a finite proper map.
  \item
    There is an algebraic curve $C\subset S$ with $C\simeq\C$ and
    a holomorphic automorphism $\phi$ of $S$ such that
    $\phi(D)\subset C$.
  \item
    There is an algebraic action of the additive group $(\C,+)$ on $S$
    and a holomorphic automorphism $\phi$ of $S$ such that
    $\phi(D)$ is contained in one orbit.
  \end{enumerate}
\end{theorem}

\begin{proof}
  $(ii)\implies(i)$: Obvious.

  $(i)\implies (ii)$: Corollary~\ref{dan-weak-strong}.

  $(i)\implies (iii)$: Proposition~\ref{tame-proper-proj}.

  $(v)\implies(iv)$: Obvious.

  $(iv)\implies (i)$: Proposition~\ref{flex-tame}.

  $(i)\implies (v)$: Let $C$ be such an orbit. It contains a tame
  discrete subset $D_0$ due to Proposition~\ref{flex-tame}.
  Thanks to Corollary~\ref{dan-tame-eq} there is a holomorphic
  automorphism $\phi$ of $X$ with $\phi(D)=D_0$ for every
  (weakly) tame discrete subset $D$ of $X$.
\end{proof}

\begin{lemma}\label{prescribe-p2}
  Let $X$ be a Danielewski surface, $D$ a discrete subset, which
  $p_1:X\to\C$ maps injectively onto a discrete subset of $\C$.

  Let $h:x\mapsto h(x)$ be an arbitrary map from $D$ to $\C$.

  Then there exists a holomorphic automorphism $\phi$ of $X$
  such that
  \begin{enumerate}
  \item
    $\pi_1\circ\phi=p_1$.
  \item
    $p_2(x)=h(x)\ \forall x\in D$.
  \end{enumerate}
\end{lemma}

\begin{proof}
  Since $p_2$ maps the $\mu_1$-orbits surjectively onto $\C$,
  for every $x\in D$ there is a complex number $\xi_x$
  such that
  \[
  p_2\left( \mu_1(\xi_x,x)\right)= h(x)
  \]
  Let $f:\C\to\C$ be a holomorphic function such that
  $f(p_1(x))=\xi_x\ \forall x\in D$.
  Now
  \[
  \phi:x\mapsto \mu_1(f(p_1(x)),x)
  \]
  is the desired automorphism.
\end{proof}

\begin{proposition}\label{dan-inj}
  Let $X$ be a Danielewski surface and let $D,\tilde D$
  be weakly tame discrete subsets
  and let $\zeta:D\to\tilde D$ be an injective map.

  Then there exists a holomorphic automorphism $\alpha$ of $X$ such that
  $\zeta(x)=\alpha(x)\ \forall x\in D$.
\end{proposition}

\begin{proof}
  Let $\phi_1,\psi_1\in\Aut(X)$ such that $p_1\circ\phi_1$
  resp.~$p_1\circ\psi_1$ maps $D$ resp.~$\tilde D$ injectively onto
  a discrete subset of $X$.
  We choose maps $h:D\to\C\setminus\{0\}$,
  $\tilde h:\tilde D\to\C\setminus\{0\}$
  such that
  \[
  \tilde h(\zeta(x))=h(x)\ \forall x\in D
  \]
  By lemma~\ref{prescribe-p2}
  there are holomorphic automorphisms $\phi_2,\psi_2$ of $Aut(X)$
  such that
  \[
  p_2(\phi_2(\phi_1(x)))=h(x)\ \forall x\in D,\quad
  p_2(\psi_2(\psi_1(x)))=\tilde h(x)\ \forall x\in \tilde D
  \]
  Since we required $h$ and $\tilde h$ to avoid the value zero,
  for every $x\in \phi_2(\phi_1(D))$
  and every $x\in \psi_2(\psi_1(\tilde D))$ the $p_2$-fiber
  of $p_2(x)$ consists of a single $\mu_1$-orbit.

  By construction
  \[
  p_2(\phi_2(\phi_1(x)))
  =p_2(\psi_2(\psi_1(\zeta(x))))\ \forall x \in D
  \]
  We choose a map
  $\gamma:D\to\C$
  such that
  \[
  (\mu_2)_{\gamma(x)}\left(\phi_2(\phi_1(x)))\right)
  =\psi_2(\psi_1(\zeta(x))))\ \forall x \in D
  \]
  Let $\beta:\C\to\C$ be a holomorphic function such that
  \[
  \beta(p_2(\phi_2(\phi_1(x))))=\gamma(x)
  \]
  and let $\eta$ be the time-one-map of the flow
  associated to $(\beta\circ p_2)\delta_2$.
  Then
  \[
   \eta\left(\phi_2(\phi_1(x)))\right)
  =\psi_2(\psi_1(\zeta(x))))\ \forall x \in D
  \]
  Thus
  \[
  \zeta(x)=\alpha(x)\ \forall x\in D
  \]
  for
  \[
  \alpha\eqdef \psi_1^{-1}\circ\psi_2^{-1}\circ\eta\circ\phi_2\circ\phi_1
  \]
  \end{proof}

\begin{corollary}\label{dan-weak-strong}
  In a Danielewski surface every weakly tame discrete subset
  is already strongly tame.
\end{corollary}

\begin{corollary}\label{dan-tame-eq}
  In a Danielewski surface any two (weakly) tame discrete subsets
  are equivalent under the action of the automorphism group, i.e.,
  given tame discrete subsets $D,D'$ there is an automorphism $\phi$
  of the surface with $\phi(D)=D'$.
\end{corollary}

\begin{theorem}\label{dan-is-RR}
  Every Danielewski surface
  \[
  X=\{(x,y,z):xy=P(z)\}
  \]
  is an $RR$-space.
\end{theorem}

\begin{proof}
  Let $\mathcal T$ be the set of weakly tame discrete subsets of $X$.
  Axiom $(RR3)$ is due to the definition of ``weakly tame''.
  $(RR4)$ is a consequence of standard
  properties of Danielewski surfaces, e.g., take
  \[
  v=P'(z)\partial_x+y\partial_z, f=y
  \]
  $(RR2)$ follows from \ref{dan-thres} and $(RR1)$ from
  \ref{dan-inj}.
  Thus $X$ satisfies all the axioms of an $RR$-space.
\end{proof}

\section{Union of tame sets}

\begin{lemma}\label{D-split}
  Let $X, Y_1,Y_2$ be manifolds, $\pi=(\pi_1,\pi_2):X\to Y_1\times Y_2$
  a proper continuous map and $D\subset X$ a discrete subset.

  Then there exists discrete subsets $D_1,D_2$ in $X$ such that
  \begin{enumerate}
  \item
    $D=D_1\cup D_2$.
  \item
    For each $i\in\{1,2\}$ the map $\pi_i$ restricts to a proper map
    from $D_i$ onto a discrete subset of $Y_i$.
  \end{enumerate}
\end{lemma}

\begin{proof}
  We choose some exhaustion functions $\rho_i$ on $Y_i$ and define
  \begin{align*}
  D_1&=\{x\in D: \rho_1(\pi_1(x))\ge\rho_2(\pi_2(x))\},\quad\\
  D_2&=\{x\in D: \rho_1(\pi_1(x))\le\rho_2(\pi_2(x))\}\\
  \end{align*}
  Since $\pi:X\to Y_1\times Y_2$ is proper,
  \[
  \rho(x)\stackrel{def}=
  \rho_1(\pi_1(x))+\rho_2(\pi_2(x))
  \]
  defines an exhaustion function on $X$.
  Hence for every $R>0$ the set $\{x\in D:\rho(x)\le R$ is finite.
  Now
  \[
  \forall x\in D_1: \rho(x)\le 2 \rho_1(\pi_1(x))
  \]
  Thus $ \rho_1(\pi_1(x))\le \frac 12R\implies \rho(x)\le R$ for $x\in D_1$.
  Hence
  \begin{align*}
  \left\{x\in D_1: \rho_1(\pi_1(x))\le \frac 12R\right\}
&  \subset
  \left\{x\in D_1: \rho(x)\le R\right\}\\
&  \subset
  \left\{x\in D: \rho(x)\le R\right\}
  \end{align*}
  is finite for all $R$.
  It follows that $\pi_1$ restricts to a proper map from $D_1$ to $Y_1$.
  In particular, $\pi_1(D_1)$ is discrete in $Y_1$.
  
  Similarily for $D_2$.
\end{proof}

\begin{theorem}\label{union}
  Let $X,Y_1,Y_2$ be Stein manifolds,
  $G_i$ (for $i\in\{1,2\}$)
  connected complex Lie groups acting on $X$,
  and $\pi=(\pi_1,\pi_2):X\to Y_1\times Y_2$ 
  be a holomorphic map such that
  \begin{enumerate}
  \item
    $\pi$ is proper, i.e., $\pi^{-1}(K)$ is compact
    for every compact subset $K\subset Y_1\times Y_2$.
  \item
    $G_i$ acts on $X$ only along the $\pi_i$-fibers.
  \item
    The $G_i$-action is not trivial.
  \item
    The automorphism group $\Aut(X)$ of $X$ acts transitively
    on $X$.
  \end{enumerate}

  Then every infinite discrete subset $D\subset X$ is the
  union of two weakly tame discrete subsets.
\end{theorem}

\begin{proof}
  Fix an infinite discrete subset $D\subset X$.
  For $i\in\{1,2\}$, let $Z_i$ be the
  fixed point set for the $G_i$-action on $X$.
  Let $Z=Z_1\cup Z_2$.
  Due to property $(iv)$ and Proposition~\ref{baire-1}
  there is an automorphism $\phi$
  of $X$ with $\phi(D)\cap Z=\{\}$.
  Thus there is no loss in generality in assuming
  that $D$ does not intersect $Z$.
  We fix non-constant holomorphic function $h_i$ on $Y_i$.
  Using Proposition~\ref{baire-1a} we may assume that both maps
  $h_i\circ\pi_i|_D:D\to\C$ are injective. As a consequence,
  the maps $\pi_i|_D:D\to Y_i$ are injective.
  
  By Lemma~\ref{unbounded-orbits} it follows that for both $i=1$ and
  $i=2$ and every $p\in D$ the $G_i$-orbit through $p$ is not
  relatively compact in $X$.
  Using Lemma~\ref{D-split} we have $D=D_1\cup D_2$
  such that $\pi_i|_{D_i}:D_i\to Y_i$ is proper.
  Thus $\pi_i|_{D_i}:D_i\to Y_i$ is proper and injective for $i=1,2$.

  Next we choose one-parameter subgroups for each $G_i$ corresponding
  to complete holomorphic vector fields $v_i$ and observe that
  (due to $(ii)$) every holomorphic function on $Y_1$ yields
  a $G_1$-invariant function on $X$ and similariy for $i=2$.

  We may thus invoke
  Corollary~\ref{crit-tame} and conclude that both $D_1$ and $D_2$
  are weakly tame.
  \end{proof}

\begin{corollary}\label{dan-union-tame}
  Let
  \[
  X=\{(x,y,z)\in\C:xy=P(z)\}
  \]
  be a Danielewski surface.

  Then every infinite discrete subset of $X$ is the union of two
  {\em tame} discrete subsets.
\end{corollary}

\begin{proof}
  Since $\Aut(X)$ acts transitively on $X$, the assertion follows
  from Proposition~\ref{dan-action} in combination with
  Theorem~\ref{union}.
  \end{proof}

\section{Appendix}

We would like to emphasize that the notion of a
{\em spreading family} as formulated in Definition~\ref{def-spread}
is a very delicate one.

For this purpose we consider the following ``toy model'' of a
family of projections.

\begin{proposition}
  Let $f$ be a non-constant holomorphic function.

  Let $\pi_t:\C^ 2\to\C$ be the family of maps parametrized by $t\in\C$
  which is defined as
  \[
  \pi_t:(x,y)\mapsto x+tf(y)
  \]

  Then $(\pi_t)$ is a {\em spreading family} for the
  Gauss measure on the parameter space $\C$
  if and  only if $f$ is a polynomial.
\end{proposition}

\begin{proof}
  $(i)$. Assume that $f$ is a polynomial.
  
  Given $r>0$, $\epsilon>0$, we choose $R>1$ such that
  the following conditions are satisfied:
  \begin{enumerate}
  \item
    $\frac{\pi r^ 2}{R}<\epsilon$.
  \item
    The Gauss measure of
    \[
    \left\{t\in\C:|t|>\frac{R-r}{\sqrt R}\right\}
    \]
    is less than $\epsilon$.
  \end{enumerate}
  If $|f(y)|\ge\sqrt R$, then
  \[
  \left |x+tf(y) \right|<r\ \implies\
  \left|\frac{x}{f(y)}+t\right|<\frac{r}{|p(y)|}\le\frac{r}{\sqrt R}
   \]
   so $t$ is contained inside a
   circle of radius less than $\frac{r}{\sqrt R}$, i.e.,
   in a disc of area $\pi\left(\frac{r}{\sqrt R}\right)^ 2$.

   If $|x|\ge R$ and $|f(y)|\le\sqrt R$, then
   \begin{align*}
   \left |x+tf(y) \right|<r
    &\implies\ 
   |tf(y)|\ge |x|-r\ge R-r\\
   &\implies\
   |t| \ge \frac{R-r}{|f(y)|]}\ge\frac{R-r}{\sqrt R}\\
   \end{align*}

   Thus $(\pi_t)$ is a spreading family provided $f$ is a polynomial.

   $(ii)$.

  Assume that $f$ is {\em not} a polynomial.
  It follows that $f$ has an essential singularity at $\infty$.
  Hence there is a sequence $z_n\to\infty$ with
  $\lim f(z_n)=0$.

  In particular, for every $r,R>0$ there exist $v=(x,y)\in\C^2$ with
  $||v||\ge|y|\ge R$ and $|x|,|f(y)|<r/2$.
  However, for such $(x,y)$ we have
  \[
  \Delta_1\subset\{ t : ||x+tf(y)||<r\}=M
  \]
  which implies that the measure of $M$ is bounded from below
  by the measure of the unit disc. Since this hold for {\em every}
  $R>0$, the family $(\pi_t)$ can not be spreading.
\end{proof}
  
\nocite{MR2126216}

\bibliography{tame}
\bibliographystyle{alpha}

\end{document}